\documentclass[12pt]{article}
\usepackage{amsfonts}
\usepackage{full page,
arcs, amssymb, amscd, amsmath, graphicx, 
wrapfig}
\newtheorem{thm}{Theorem}[section]
\newtheorem{defn}[thm]{Definition}
\newtheorem{prop}[thm]{Proposition}
\newtheorem{cor}[thm]{Corollary}

\newtheorem{rema}[thm]{Remark}

\newtheorem{prob}[thm]{Problem}

\newtheorem{conj}[thm]{Conjecture}

\newcommand{\nn}{\nonumber \\}

 \newcommand{\res}{\mbox{\rm Res}}

\renewcommand{\hom}{\mbox{\rm Hom}}

\newcommand{\lbar}{\bigg\vert}
\newcommand{\tr}{\mbox{\rm Tr}}

\newcommand{\Y}{\mathcal{Y}}
\newcommand{\C}{\mathbb{C}}
\newcommand{\Z}{\mathbb{Z}}
\newcommand{\R}{\mathbb{R}}

\newcommand{\N}{\mathbb{N}}

\newcommand{\one}{\mathbf{1}}

\title{ {\bf Representation theory of vertex operator algebras and 
orbifold conformal field theory} }
\date{}
\author{Yi-Zhi Huang}

\begin{document}

\bibliographystyle{alpha}
\maketitle
\begin{abstract}
We discuss some  basic problems and conjectures in a program to construct 
general orbifold conformal field theories using the representation theory of 
vertex operator algebras. We first review a program to construct conformal field theories.
We also clarify some misunderstandings on vertex operator algebras, 
modular functors and intertwining operator algebras.
Then we discuss some basic open problems and conjectures in 
mathematical orbifold conformal field theory. Generalized twisted modules and their variants, 
their constructions and some existence results are reviewed. 
Twisted intertwining operators and their basic properties are also reviewed. 
The conjectural properties in the basic open problems and conjectures mentioned above
are then formulated precisely and explicitly. Some thoughts of the author 
on further developments of orbifold conformal field theory are also discussed.
\end{abstract}

\renewcommand{\theequation}{\thesection.\arabic{equation}}
\renewcommand{\thethm}{\thesection.\arabic{thm}}
\setcounter{equation}{0}
\setcounter{thm}{0}
\section{Introduction}

Two-dimensional orbifold conformal field theories are two-dimensional 
conformal field theories constructed from 
known conformal field theories and their automorphisms. The first example of 
two-dimensional orbifold conformal field theories, 
the moonshine module, was  constructed by 
Frenkel, Lepowsky and Meurman \cite{FLM1} \cite{FLM2} \cite{FLM} in mathematics. 
In this construction, twisted vertex operators studied by Lepowsky \cite{L} 
played an important role. The systematic study of two-dimensional orbifold conformal field theories 
in string theory was started by Dixon, Harvey, Vafa and Witten  in  \cite{DHVW1} and \cite{DHVW2}. 
(Since we do not discuss higher-dimensional conformal field theories in this paper, 
as usual we shall omit the word ``two-dimensional" unless we want to emphasize the 
special feature of two-dimensional theories.)
Since then, orbifold conformal field theory has been developed in mathematics as 
an important part of mathematical conformal field theory. 

Orbifold conformal field theory is not just 
a mathematical procedure to obtain new examples of conformal field theories. 
More importantly, we expect that it will provide a powerful approach to solve mathematical 
problems and prove mathematical conjectures after it is fully developed. 
For example, one of the most important conjecture in the theory of vertex operator 
algebra and mathematical conformal field theory is the uniqueness of the moonshine module formulated
by Frenkel, Lepowsky and Meurman \cite{FLM}. The moonshine module is constructed 
as an orbifold conformal field theory from the Leech lattice vertex operator algebra and 
its automorphism induced from the point reflection in the origin 
 of the Leech lattice (see  \cite{FLM}). From the 
viewpoint of orbifold conformal field theory, the uniqueness means that every vertex operator algebra
obtained as an orbifold conformal field theory satisfying the three conditions in the uniqueness
conjecture must be isomorphic to the moonshine module as a vertex operator algebra. 
In particular, we have to  study general orbifold conformal field theories
satisfying the three conditions in this uniqueness conjecture. Such a study in turn means that 
we have to develop a general orbifold conformal field theory.

We would like to emphasize that to solve mathematical problems and prove mathematical 
conjectures, it is important to develop mathematical 
techniques that might not be most interesting in physics. For example,
existence of orbifold conformal field theories must be proved, not assumed. Since the
state spaces of orbifold conformal field theories are always infinite dimensional, it is important to prove, not assume,
the convergence of series appearing in mathematical constructions and proofs. 

On the other hand, the ideas from physics will undoubtedly play a crucial role in this development of general
orbifold conformal field theory. In the last thirty years, mathematicians have witnessed the power of  ideas 
coming from 
quantum field theory. Quantum field theory not only provides beautiful interpretations 
of many deep mathematical results, but also suggests  new approaches to 
many difficult mathematical problems. Since topological quantum field theory has been 
understood quite well in mathematics, the approaches developed using 
topological quantum field theory have been very successful in solving mathematical problems. 
On the other hand, although the ideas and conjectures from nontopological quantum field theories
such as two-dimensional conformal field theory and four-dimensional Yang-Mills theory have led to solutions of some 
important problems in mathematics, the nontopological quantum-field-theoretic 
approaches themselves have not been completely developed and some related mathematical 
problems were solved using approaches different from the still-to-be-developed 
quantum-field-theoretic approaches. In the author's opinion, to understand completely 
the beautiful mathematical structures suggested by nontopological quantum field theories and 
to solve completely the deep mathematical problems coming from nontopological quantum-field-theoretic
ideas, it is necessary to construct these nontopological quantum field theories mathematically
and to develop them into new mathematical approaches and tools for mathematical problems. 

Conformal field theory is the main nontopological quantum field theory that
has been developed substantially in mathematics. We expect that after the mathematical constructions
are completed, conformal field theory can be further developed 
to have the power to solve a number of difficult 
mathematical problems. Orbifold conformal field theory as an important part of 
conformal field theory
will be crucial in this further development. 

In the construction of conformal field theories 
using the representation theory of vertex operator algebras, chiral conformal field theories are in fact 
the same as the theory of intertwining operators. Intertwining operators among modules 
for a vertex operator algebra together with their properties can be used as a working definition of 
chiral conformal field theory.  See Section 2 for a discussion about a program to construct conformal 
field theories using the representation theory of vertex operator algebras. This program is very successful 
for rational conformal field theories and their logarithmic generalizations. 

We know much less about orbifold conformal field theories.
In fact, only orbifold conformal field theories associated to finite cyclic groups 
have been fully constructed by van Ekeren, M\"{o}ller
and Scheithauer \cite{EMS} and by M\"{o}ller \cite{M}
based on the basic constructions and results on rational conformal field theories 
obtained  by the author in \cite{H-diff-eqn}, \cite{H-modular}, \cite{H-verlinde-conj} and \cite{H-rigidity}. 
Although there are  examples of nonabelian orbifold conformal field theories (see for example
the one constructed from a lattice of rank $72$  and a nonabelian finite 
group of its automorphisms   by Gem\"{u}nden and Keller \cite{GK}),
the full development of mathematical orbifold conformal theory,
especially  those associated to nonabelian automorphism groups, 
is still in the beginning stage. 

In the present paper, we discuss some  basic problems and conjectures in a program to construct 
general orbifold conformal field theories using the representation theory of 
vertex operator algebras. Most of the problems and conjectures discussed in this paper
were proposed but were not formulated explicitly in the paper \cite{H-problems}. 
These problems and conjectures are the foundation for future mathematical developments of 
orbifold conformal field theory. As in the case of conformal field theories mentioned above, 
our approach is to view chiral orbifold conformal field theories as the theory of twisted intertwining 
operators. We use twisted intertwining operators among twisted modules for a vertex operator algebra
and a group of its automorphisms and their (conjectural) properties as a working definition of chiral 
orbifold conformal field theory. Constructing a chiral orbifold conformal field theory is the same as
proving the conjectural properties of twisted intertwining operators. 

At the end of Section 2, we also clarify some misunderstandings on vertex operator algebras, 
modular functors and intertwining operator algebras, which for many years,
in the author's opinion,  have often become an obstruction to
 the mathematical development of conformal field theory 
and the representation theory of vertex operator algebras.

The present paper is organized as follows: In Section 2, we review a program to construct 
conformal field theories using the representation theory of vertex operator algebras. 
The basic  problems and conjectures on a general construction of orbifold conformal field theories
are given in Section 3. These problems and conjectures are formulated using 
suitable twisted modules and twisted intertwining operators. So in Sections 4 and 5, we 
recall the basic notions and general results about suitable twisted modules and twisted intertwining operators, 
respectively. In Section 6, we formulate and discuss in details the conjectural properties of twisted intertwining operators
stated in the problems and conjectures in Section 3. Although there is no result in Section 6 and the 
formulations are simple generalizations of those in the untwisted case, 
the material in this section has not appeared in any other 
publications. In Section 7, we discuss some thoughts of the author 
on further developments based on these problems and conjectures.

\paragraph{Acknowledgments}
The author is grateful to Drazen Adamovic and Antun Milas for their invitation to the conference 
``Representation Theory XVI" held in Dubrovnik in June 2019. 

\renewcommand{\theequation}{\thesection.\arabic{equation}}
\renewcommand{\thethm}{\thesection.\arabic{thm}}
\setcounter{equation}{0}
\setcounter{thm}{0}
\section{A program to construct  conformal field theories}

Vertex operator algebras were introduced in mathematics by Borcherds \cite{B} and 
Frenkel, Lepowsky and Meurman \cite{FLM} in the study of the moonshine module
underlying the Monster finite simple group. They
are also algebraic structures appearing  in physics  in the study of 
conformal field theory in the work of Belavin, Polyakov and Zamolodchikov \cite{BPZ} and 
Moore and Seiberg \cite{MS}. Mathematically, Kontsevich and G. Segal \cite{S}
formulated a definition of conformal field theory.  
A conformal field theory is,  roughly speaking,  an algebra 
over the operad of Riemann surfaces 
with parametrized boundaries satisfying certain additional conditions 
(see \cite{S},  \cite{H-program-1}, \cite{H-book}
and \cite{H-mod-f-int-alg}). Conformal 
field theories can be constructed and 
studied using  intertwining operators (or chiral vertex operators in \cite{MS}) among modules 
for vertex operator algebras satisfying suitable conditions. 

The first main problem in the mathematical study of conformal field theory is to give a mathematical construction. 
A program to construct conformal field theories using 
the representation theory of vertex operator
algebras has been quite successful in the past thirty years. 
In this program,  the construction 
of a chiral conformal field theory can be divided into the following steps: (i) Construct a 
vertex operator algebra and study their modules. (ii) Prove the associativity 
 (or operator product expansion) of intertwining operators among modules
for this vertex operator algebra. (iii) Prove the modular 
invariance of products of intertwining operators. (iv) Prove a higher-genus convergence property
(this  is still a conjecture even for a rational chiral conformal field theory). 
To construct a full conformal field theory, one needs to further construct 
a nondegenerate bilinear form or Hermitian form on the space of intertwining operators
satisfying certain conditions so that 
a chiral conformal field theory and an antichiral conformal field theory can be put together to 
give a full conformal field theory.  See \cite{H-program-1}, 
\cite{H-book}, \cite{H-program-2} and \cite{H-blog} for expositions on 
this program and \cite{H-problems} for some of the remaining main open problems
in this program. 

Many results in the representation theory of vertex operator 
algebras also hold for more general grading-restricted vertex algebras
and M\"{o}bius vertex algebras. So here we recall the 
the notion of vertex operator algebra by first recalling these notions. 
See  \cite{FHL}, and \cite{HLZ1}, \cite{H-2-const} and \cite{H-exist-twisted-mod} 
for details. 
Roughly speaking, a grading-restricted vertex algebra is a $\Z$-graded vector space $V=\coprod_{n\in \Z}V_{(n)}$
satisfying the grading-restriction conditions $\dim V_{(n)}<\infty$ for $n\in \Z$ and 
$V_{(n)}=0$ for $n$ sufficiently negative, equipped with a vertex operator map
\begin{eqnarray*}
Y_{V}: V\otimes V&\to& V[[z, z^{-1}]],\nn
u\otimes v&\mapsto& Y(u, z)v.
\end{eqnarray*}
(analogous to the multiplication for an associative algebra but parametrized by $z$), a vacuum $\one\in V_{(0)}$ 
(analogous to the identity in an associative algebra)  satisfying axioms 
similar to those for commutative associative algebras and also axioms 
on the derivatives of vertex operators with respect to $z$
and the meromorphicity of the correlation functions obtained from 
products of vertex operators. A M\"{o}bius vertex algebra is a grading-restricted vertex 
algebra $V$ equipped with a compatible structure of $\mathfrak{sl}_{2}$-module.
A vertex operator algebra is a M\"{o}bius vertex algebra
 equipped with a compatible conformal vector such that the components of the vertex operator 
of the conformal vector satisfy the Virasoro relations and some other properties. 
These algebras are in fact analogues of commutative associative algebras
with their multiplications controlled by Riemann spheres with three punctures and local 
coordinates vanishing at the punctures. See \cite{H-book} for details. 

We now discuss axioms for these algebras. 
For example, we have two  axioms for the vacuum: For $u\in V$, 
$Y(\mathbf{1}, z)u=u$ and $\lim_{z\to 0}Y(u, z)\mathbf{1}=u$. These are analogues 
of the axioms for the identity in an associative algebra.
The main axiom is the associativity (the analogue of the associativity for an associative algebra). 
It says that  for $u_{1}, u_{2}, u_{3}\in V$ and $u'\in V'$ (the graded dual of $V$), 
$$\langle u', Y(u_{1}, z_{1})Y(u_{2}, z_{2})u_{3}\rangle$$ and 
$$\langle u', Y(Y(u_{1}, z_{1}-z_{2})u_{2}, z_{2})u_{3}\rangle$$
are absolutely convergent in the regions $|z_{1}|>|z_{2}|>0$ and 
$|z_{2}|>|z_{1}-z_{2}|>0$, respectively, to a common rational function in $z_{1}, z_{2}$ with the only possible poles at
$z_{1}=0$, $z_{2}=0$ and $z_{1}=z_{2}$.
Another main axiom is the commutativity (the analogue of the commutativity 
for a commutative associative algebra in a more subtle sense). 
Roughly speaking, it says that for $u_{1}, u_{2}\in V$, we require that the rational functions
obtained by analytically extending 
$$\langle u', Y(u_{1}, z_{1})Y(u_{2}, z_{2})u_{3}\rangle$$ and 
$$\langle u', Y(u_{2}, z_{2})Y(u_{1}, z_{1})u_{3}\rangle$$
are the same. In addition, we have the $L_{V}(-1)$-derivative property 
$$\frac{d}{dz}Y(u, z)=Y(L_{V}(-1)u, z)$$
for $v\in V$, where 
$$L_{V}(-1)v=\lim_{x\to 0}Y_{V}(v, z)\one$$
for $v\in V$.This property means that $L_{V}(-1)$ corresponding 
to the derivative or infinitesimal translation of the variable $z$) 
 For a M\"{o}bius vertex algebra, let $L_{V}(1), L_{V}(0), L_{V}(-1)$ be the  operators  giving
the $\mathfrak{sl}_{2}$-module structure. Then 
$$L_{V}(-1)v=\lim_{z\to 0}Y_{V}(v, z)\one$$
for $v\in V$, $L_{V}(0)$ is the operator giving the grading of $V$ and 
$$[L_{V}(1), Y_{V}(v, z)]=Y_{V}(L_{V}(1)v, z)+2zY_{V}(L_{V}(0)v, z)
+z^{2}Y_{V}(L_{V}(-1)v, z).$$
For a vertex operator algebra, we have the axioms for the conformal symmetry.  For example, if we let $L_{V}(n)
=\res_{z} z^{n+1}Y_{V}(\omega, z)$, then 
$$[L_{V}(m), L_{V}(n)]=(m-n)L_{V}(m+n)+\frac{c}{12}(m^{3}-m)\delta_{m+b, 0}$$
(the Virasoro relation). Also the Virasoro operators $L_{V}(1)$, $L_{V}(0)$ and $L_{V}(-1)$ should be the same 
as those operators for the underlying M\"{o}bius vertex algebra. 

A module for a vertex operator algebra $V$ is roughly speaking a $\C$-graded
vector space and a vertex operator map 
$Y_{W}: V\otimes W\to W[[z, z^{-1}]]$
satisfying all the axioms that still make sense. 
For three modules $W_{1}$, $W_{2}$ and $W_{3}$, an intertwining operator of type $\binom{W_{3}}{W_{1}W_{2}}$
is a linear map $\Y: W_{1}\otimes W_{2}\to W_{3}\{z\}[\log z]$ 
(here $W_{3}\{z\}$ means series in complex powers of $z$  with coefficients in $W_{3}$) 
satisfying 
all axioms for modules that still make sense. 

The intertwining operators of type $\binom{W_{3}}{W_{1}W_{2}}$ form a vector space. 
This is in fact the space of conformal blocks on the Riemann sphere with three marked points 
labeled with the equivalence classes of the modules $W_{1}$, $W_{2}$ and $W_{3}'$. 
Its dimension is called a fusion rule. Intertwining operators were called chiral vertex operators 
in \cite{MS} and were introduced mathematically in \cite{FHL}. 
If a set of modules for a vertex operator algebra equipped with 
subspaces of intertwining operators 
among modules in this set satisfying the associativity 
and commutativity, we call it
an intertwining operator algebra (see \cite{H-program-1}, \cite{H-mod-f-int-alg} and 
\cite{H-generalizedrationality} and see also \cite{DL} for a special type of intertwining operator algebras
called abelian intertwining algebras).  Here by associativity of intertwining operators, 
we mean, roughly, for any intertwining operators 
$\Y_{1}$ and $\Y_{2}$, there exist intertwining operators 
$\Y_{3}$ and $\Y_{4}$ such that for $w_{1}, w_{2}$, $w_{3}$ and $w_{4}'$ in suitable modules,
$$\langle w_{4}', \Y_{1}(w_{1}, z_{1})\Y_{2}(w_{2}, z_{2})w_{3}\rangle$$
and 
$$\langle w_{4}', \Y_{3}(\Y_{4}(w_{1}, z_{1}-z_{2})w_{2}, z_{2})w_{3}\rangle$$
are absolutely convergent in the regions $|z_{1}|>|z_{2}|>0$ and 
$|z_{2}|>|z_{1}-z_{2}|>0$, respectively.
Moreover, these functions can be analytic extended to a common
multivalued analytic function with the only possible singular points at 
$z_{1}=0, z_{2}=0$ and $z_{1}=z_{2}$. The associativity of intertwining operators 
is equivalent to the operator product expansion of chiral vertex operators,
one of the two major assumptions or conjectures in the important work \cite{MS} of Moore and Seiberg.
Mathematically it was first introduced and proved under suitable conditions in \cite{tensor4}. 
By commutativity of intertwining operators, we mean, 
roughly, for any intertwining operators 
$\Y_{1}$ and $\Y_{2}$, there exist intertwining operators 
$\Y_{3}$ and $\Y_{4}$ such that 
$$\langle w_{4}', \Y_{1}(w_{1}, z_{1})\Y_{2}(w_{2}, z_{2})w_{3}\rangle$$
and 
$$\langle w_{4}', \Y_{3}(w_{2}, z_{2})\Y_{4}(w_{1}, z_{1})w_{3}\rangle$$
are analytic extensions of each other. The commutativity of intertwining operators 
is an easy consequence of the associativity and skew-symmetry of intertwining operators
(see \cite{H-generalizedrationality}).

Intertwining operator algebras give vertex tensor category structures which in turn give
braided tensor category structures. 
Vertex tensor categories can be viewed as analogues of
symmetric tensor categories with their tensor product bifunctors controlled 
by Riemann surfaces with three punctures and local coordinates vanishing at the punctures. 
See \cite{HL-vertex-tensor-cat}, \cite{H-rigidity} and \cite{HLZ8} for details.

Intertwining operator algebras are equivalent to chiral genus-zero conformal field theories (see \cite{H-mod-f-int-alg}). 
In the program to construct conformal field theories using the representation theory of 
vertex operator algebras, after the first step of constructing a vertex operator algebra and studying 
its modules are finished,  the second step of proving the associativity (or operator product expansion) 
of intertwining operators is equivalent to constructing an intertwining operator algebra.
To prove the associativity, the main properties that need to be established first are
a convergence and extension property for products of an 
arbitrary number of intertwining operators (see \cite{H-diff-eqn} in the rational case, \cite{HLZ7} 
for the adaption of the proof in \cite{H-diff-eqn} in the logarithmic case and \cite{Y} for a
generalization to vertex algebras with infinite-dimensional homogeneous subspaces with respect to
 conformal weights but with finite-dimensional homogeneous subspaces 
with respect to an additional horizontal grading)
and a property stating that suitable lower-bounded generalized modules are  in the category of 
modules that we start with (see \cite{tensor4} for the rational case, \cite{HLZ7} and \cite{H-zhang-correction}
for the generalization to the logarithmic case). 
The associativity is proved using these properties (see \cite{tensor4} for the rational case and \cite{HLZ6}
for the generalization to the logarithmic case). 

The next step is to construct the chiral genus-one conformal field theories. 
The main properties that need to be established are the convergence and analytic extensions of 
$q$-traces or pseudo-$q$-traces of products of geometrically-modified 
intertwining operators (see \cite{H-modular} for the rational case and \cite{F1}
and  \cite{F}
for the generalization to the logarithmic case) 
and the modular invariance of  the analytic extensions of  these $q$-traces 
or pseudo-$q$-traces (see \cite{H-modular}
for the rational case).
The genus-one associativity and 
commutativity are easy consequences of the convergence and analytic extensions of these $q$-traces
or pseudo-$q$-traces
(see \cite{H-modular}  and \cite{F1} and \cite{F}
for the generalization to the logarithmic case). 

The main open problem in the construction of chiral higher-genus rational conformal 
field theories is the convergence of multi $q$-traces of  products of geometrically-modified intertwining operators.
The invariance under the action of the mapping class groups is an easy consequence 
of this convergence, the associativity of intertwining operators and the modular invariance 
of the $q$-traces of products of geometrically-modified  intertwining operators. The future solution of this problem 
will depend on the further study of the moduli space of Riemann surfaces with parametrized boundaries, including 
in particular the study of a conjecture by the author on  meromorphic functions on this moduli space.
See \cite{RS1}--\cite{RS3} and \cite{RSS1}--\cite{RSS5} for  results on this moduli space. 

Another problem is the construction of locally convex topological completions of modules for the vertex operator algebra
such that intertwining operators and higher-genus correlation functions give maps between these 
completions. Such completions of vertex operator algebras and their modules using only the correlation functions
given by the algebras and modules were given in \cite{H-completion-1} and \cite{H-completion-2}. 
If we assume the convergence of multi $q$-traces of  products of geometrically-modified intertwining operators
discussed above, then the same method as in \cite{H-completion-1} and \cite{H-completion-2}
works when we add those elements coming from genus-zero and genus-one correlation functions 
obtained using intertwining operators. The author conjectured that these completions 
obtained using all genus-zero and genus-one correlation functions are the same as the Hilbert 
space completions if the chiral conformal field theory is 
unitary (see \cite{H-oberwolfach}). 

The discussions above are about the construction of chiral conformal field theories.
We also need to construct full conformal field theories and open-closed conformal 
field theories. To construct a genus-zero full rational conformal field theory, the main property
that needs to be proved is the nondegeneracy of an invariant bilinear from on the space of 
intertwining operators. This nondegeneracy in fact needs a formula used by the author to prove
the Verlinde formula in \cite{H-verlinde-conj} or equivalently 
the rigidity of the braided tensor category of modules for the vertex operator algebra
proved in \cite{H-rigidity}. The construction of genus-one and higher-genus 
full conformal field theories can then be obtained easily from genus-zero 
full conformal field theories
and chrial genus-one and higher-genus conformal field theories (see \cite{HK1} and \cite{HK2}).

Finally, to construct open-closed conformal field theories, one needs to construct 
open-string vertex operator algebras from intertwining operator algebras (see \cite{HK0}). 
Then one has to prove that with the choices of the open-string vertex algebra (the open part)
and the full conformal field theory (the closed part), Cardy's condition on the compatibility between 
the closed part and open part is satisfied (see \cite{K1} and \cite{K2}). 

From the discussion above, we see that intertwining operators are the main objects to study in 
conformal field theory and also in the representation theory of vertex operator algebras. In fact,
conformal field theory is essentially the theory of intertwining operators. 
Therefore  intertwining operators with all their properties can be viewed as a 
working definition of chiral conformal field theory. Constructing a chiral conformal field theory is 
the same  as proving all the properties of intertwining operators. 
This will also be our approach in this paper to the construction of orbifold conformal field theories. 

Before we end this section,  the author would like to correct some misunderstandings about conformal field theories
and vertex operator algebras because for a long time, opinions formed based on 
these and other misunderstandings 
have often been used mistakenly
by journals, organizations and the mathematical community to evaluate researches in this area.

The first misunderstanding is about the role of vertex operator algebras 
in conformal field theory. A vertex operator algebra in general is certainly not even a chiral conformal field theory. 
This can be seen easily  from the modular invariance of the space of 
characters of integrable highest weight modules for  affine Lie algebras (see \cite{Ka})
and from Zhu's modular invariance theorem (see \cite{Z})
on $q$-traces of vertex operators acting on modules for a vertex operator algebra. Also many poweful methods 
used to study vertex operator algebras do not work for intertwining operators. For example, 
since vertex operator algebras involve only rational functions,  the method of multiplying a 
polynomial to cancel the denominator of a rational function works very well. But this method in general
does not work for products of at least two intertwining operators. Also for rational functions, 
one can use the global expressions of rational functions instead of analytic extensions but for 
multivalued functions obtained from products of at least two intertwining operators, 
one has to carefully use analytic extensions to obtain the correct results.
In fact, fatal mistakes occurred in papers
published in major mathematical journals claiming to simplify major results on intertwining operators
without using complex analysis  exactly because the methods that work only 
for vertex operator algebras were applied to the study of 
products of two intertwining operators. One mistake is
to assume that intertwining operators involve only integral powers of the variable, 
which  is not even true for non-self-dual lattice vertex operator algebras,
the simplest 
minimal model of central charge $\frac{1}{2}$ and the simplest Wess-Zumino-Witten 
models for the Lie algebra $\mathfrak{sl}_{2}$. Another mistake is to define
maps without using analytic extensions. 
When working with multivalued analytic functions (not rational functions), one has to use analytic extensions to define 
a number of maps. Without a careful use of analytic extensions, one cannot even prove that an arbitrarily defined 
map is linear, not to mention many other properties that these maps should satisfy. 

Another more subtle but also more important fact is that a vertex operator algebra in general
does not even determine a chiral conformal field theory uniquely.  Instead, it is in fact 
the choices of modules and intertwining operators that determine uniquely such a theory. 
Therefore to construct a chiral conformal field theory, though we need to start with a vertex operator algebra, 
the more crucial part is to choose a category of modules and spaces of intertwining operators. 
One simple example is the chiral conformal field theory associated to irrational tori. 
The vertex operator algebra for such a chiral conformal field theory is a Heisenberg vertex operator algebra,
which is also the vertex operator algebra for the conformal field theory associated to the corresponding Euclidean space. 
This vertex operator algebra alone does not lead us to a unique chiral conformal field theory since 
all different irrational tori and the Euclidean space in the same dimension has the same vertex operator algebra. We have 
to choose the category of modules for this Heisenberg vertex operator algebra to be the category
of finite direct sums of irreducible modules generated by 
the eigenfunctions of the Laplacian on the given irrational torus 
and the space of  intertwining operators among the irreducible modules
to be the spaces of all intertwining operators among these modules. Then we obtain the chiral 
conformal field theory associated to this particular irrational torus. 

The second misunderstanding is about modular functors and modular tensor categories. 
Modular functors are operads formed by holomorphic vector bundles over the moduli space of 
Riemann surfaces with parametrized boundaries satisfying certain additional conditions
and  chiral conformal field theories are algebras over such operads satisfying additional conditions
(see \cite{S} and also \cite{H-program-1}
and \cite{H-mod-f-int-alg} for the genus zero case). But modular functors themselves 
are not conformal field theories.
To construct a chiral conformal field theory, one has to construct a modular functor
with a nondegenrate bilinear form. But 
a modular functor with a nondegenrate bilinear form 
alone does not give a conformal field theory. Instead, a modular functor with a nondegenrate bilinear form gives a 
three-dimensional topological field theory. Similarly, modular tensor categories also 
give only three-dimensional topological field theories and is equivalent to modular functors with 
nondegenerate bilinear forms. They are far from conformal field theories.
Certainly modular functors and modular tensor categories are very useful in the study 
of conformal-field-theoretic structures. 
But this happens only when we already proved a lot of results, for example, the convergence, 
associativity (operator product expansion) and modular invariance, about intertwining operators or some 
equivalent structures. For example, the modular tensor categories for the Wess-Zumino-Witten models
can be constructed using representations of quantum groups. The fusion coefficients of 
these tensor categories are indeed given by the Verlinde formula. But these modular tensor categories do
not give us the convergence, associativity (operator product expansion)
 and modular invariance for intertwining operators. 

The third misunderstanding is about the history of intertwining operator algebras.
Original papers, books and proposals on intertwining operator algebras have often been 
rejected  based on wrong  claims that some much later or nonexistent mathematical works 
 already had this notion or some results before 
intertwining operator algebras were actually introduced and studied. 
To correct this misunderstanding, 
we give a brief history of intertwining operator 
algebras here. (To be accurate, some of the years below are the years that the papers appeared in 
the arXiv, not the years that they were published.)

In 1984, operator product expansion of chiral conformal fields was studied by Belavin, 
Polyakov and Zamolodchikov in \cite{BPZ}. In 1988, by assuming that two major conjectures---the 
operator product expansion and modular invariance (certainly including implicitly the corresponding convergence)
of chiral vertex operators (equivalent to intertwining operators in mathematics)---hold, 
Moore and Seiberg derived a set of polynomial equations 
which corresponds to a modular tensor category in the sense of Turaev \cite{T1}
and obtained the Verlinde formula \cite{V} as a consequence. In 1992, Dong and Lepowsky 
 in \cite{DL}  introduced a special type of intertwining operator algebras
called abelian intertwining algebras for which the corresponding braid group representations
are one dimensional and gave examples constructed from lattices. In 1995, 
the author in \cite{tensor4} formulated and proved the associativity of intertwining
operators assuming that a convergence and extension property for intertwining operators and another 
algebraic condition hold. 
In particular, the operator product expansion of intertwining operators was proved assuming these conditions.
At the same time in 1995, the author in \cite{H-minimal-model} proved the 
convergence and extension property and the other algebraic condition needed in 
\cite{tensor4} for minimal models. In the same year, the author 
introduced in \cite{H-program-1} the mathematical notion of intertwining operator algebra using the associativity  of 
intertwining operators and discussed its role in the 
construction of conformal field theories in the sense of Kontsevich and Segal \cite{S}. 

In 1997
the author proved a generalized rationality and a Jacobi identity for intertwining operator algebras. 
In the same year, Lepowsky and the author in \cite{HL-affine}
proved the convergence and extension property and the other 
algebraic condition needed in \cite{tensor4} for the Wess-Zumino-Witten models. 
Also in the same year, the author constructed genus-zero modular functors from 
intertwining operator algebras and proved that intertwining operator algebras are 
algebras over the partial operads of such genus-zero modular functors. 
In 1999 and 2000, Milas and the author in \cite{HM1} and \cite{HM2}
proved the convergence and extension property and the other algebraic condition needed in 
\cite{tensor4}  for the Neveu-Schwarz sectors of $N$=1 and $N=2$  superconformal minimal 
models,  respectively. In 2001, the author in \cite{H-cft-lattice}
introduced a notion of dual of an intertwining operator algebra analogous to 
the dual of a lattice and the dual of a code such that the dual of a vertex operator algebra 
satisfying suitable conditions is the intertwining operator algebra obtained from all intertwining operators
among all irreducible modules. In 2002, the author proved 
 the convergence and extension property and the other condition 
needed in \cite{tensor4} for vertex operator algebras 
for which irreducible modules are $C_{1}$-cofinite (weaker than $C_{2}$-cofnite)  and 
$\N$-gradable weak modules are completely reducible.  The convergence 
and extension property proved in this paper in fact also holds for vertex operator algebras 
for which grading-restricted generalized modules are  $C_{1}$-cofinite 
but might not be 
completely reducible (see \cite{HLZ7} for a discussion about this fact). 
In 2003, the author proved the modular invariance for intertwining operator algebras on the
direct sum of all (inequivalent representatives of) irreducible modules for a vertex operator algebra
satisfying natural finite reductivity conditions, including in particular the $C_{2}$-cofniteness 
condition and the complete reducibility of $\N$-gradable weak modules. 
A number of these results on intertwining operator algebras were generalized 
to the logarithmic and other cases by Lepowsky, Zhang and the author in \cite{HLZ0}--\cite{HLZ8}, by 
the author in \cite{H-cofiniteness} and \cite{H-zhang-correction},  by Fiordalisi in
\cite{F1} and \cite{F} and by Yang in \cite{Y}. 

\renewcommand{\theequation}{\thesection.\arabic{equation}}
\renewcommand{\thethm}{\thesection.\arabic{thm}}
\setcounter{equation}{0}
\setcounter{thm}{0}
\section{Basic  open problems and conjectures in orbifold conformal field theory}

In this section, we discuss some basic open problems and conjectures in orbifold conformal field theory. 

Roughly speaking, given a conformal field theory and a group of automorphisms of this theory, 
we would like to construct a conformal field theory that can be viewed as the original conformal field theory 
divided by this group of automorphisms. To formulate this notion precisely, 
we assume that the given conformal field theory is constructed using the representation theory of 
vertex operator algebras, as is discussed in the preceding section, and the group of automorphisms
is a group of automorphisms of the vertex operator algebra. 

Here is the main open problem on the construction of orbifold conformal field theories:

\begin{prob}\label{main-prob}
Given a vertex operator algebra $V$ and a group $G$ of automorphisms of $V$,
construct and classify all the conformal field theories whose vertex operator algebras 
contain the fixed point vertex operator algebra
$V^{G}$ as subalgebras. 
\end{prob}

It is in fact very difficult to study $V^{G}$-modules 
and intertwining operators among $V^{G}$-modules. On the other hand, twisted 
$V$-modules and twisted intertwining operators among twisted $V$-modules are analogues 
of $V$-modules and intertwining operators among $V$-modules, it is easier to study 
these than $V^{G}$-modules 
and intertwining operators among them. We expect that 
$V^{G}$-modules can all be obtained from  twisted 
$V$-modules (see Theorem 7.1 in \cite{H-exist-twisted-mod} and Theorem \ref{v-g-extension} below
for lower-bounded $V^{G}$-modules in the case that $G$ is the cyclic 
group generated by an automorphism $g$ of $V$). Thus our conjectures and problems below will be mainly on twisted intertwining 
operators among twisted modules. 

Now we state the first main conjecture on orbifold conformal field theories (see \cite{H-problems}):

\begin{conj}\label{twisted-int-op}
Assume that $V$ is a simple vertex operator algebra 
satisfying the following conditions:
\begin{enumerate}
\item $V_{(0)}=\C\one$, $V_{(n)}=0$ 
for $n<0$ and the contragredient $V'$, as a $V$-module, is equivalent to $V$.

\item $V$ is $C_{2}$-cofinite, that is, $\dim V/C_{2}(V)<\infty$, where $C_{2}(V)$ is the subspace of 
$V$ spanned by the elements of the form $\res_{x}x^{-2}Y(u, x)v$ for $u, v\in V$ and $Y: V\otimes V\to V[[x, x^{-1}]]$
is the vertex operator map for $V$.

\item Every grading-restricted generalized $V$-module is completely reducible.
\end{enumerate}
 Let $G$ be a finite group 
of automorphisms of $V$. Then the twisted intertwining operators among the $g$-twisted $V$-modules for all $g\in G$
satisfy the associativity, commutativity and modular invariance property.
\end{conj}

The following conjecture (see also \cite{H-problems})
 is a  consequence of Conjecture \ref{twisted-int-op} (cf. Example 5.5 in \cite{Ki}):

\begin{conj}\label{g-crossed}
Let $V$ be a vertex operator satisfying the three conditions in Conjecture \ref{twisted-int-op}
and let $G$ be a finite group of automorphisms of $V$. The the category of $g$-twisted $V$-modules for 
all $g\in G$ is a $G$-crossed modular tensor category (see Turaev \cite{T} for $G$-crossed tensor categories). 
\end{conj}

These two conjectures given in \cite{H-problems} are both under the complete reducibility assumption and are also about 
finite groups of automorphisms of $V$. In the case that grading-restricted generalized $V$-modules 
are not complete reducible or $G$ is not finite, we have the following conjectures and problems:

\begin{conj}
Let $V$ be a vertex operator algebra satisfying the first two conditions in Conjecture \ref{twisted-int-op}
and let $G$ be a finite group of automorphisms of $V$. Then the twisted  intertwining operators 
among the grading-restricted generalized $g$-twisted $V$-modules for all $g\in G$
satisfy the associativity, commutativity and modular invariance properties.
\end{conj}

\begin{conj}
Let $V$ be a vertex operator algebra satisfying the first two conditions in Conjecture \ref{twisted-int-op}
and let $G$ be a finite group of automorphisms of $V$. 
Then the category of grading-restricted generalized $g$-twisted $V$-modules for 
all $g\in G$ has a natural structure of $G$-crossed 
tensor category satisfying additional properties. 
\end{conj}

\begin{prob}
Let $V$ be  a vertex operator algebra and let $G$ be a group of automorphisms of $V$. If $G$ is an infinite group,
under what conditions do the twisted intertwining operators 
among the grading-restricted generalized $g$-twisted $V$-modules for all $g\in G$
satisfy the associativity, commutativity and modular invariance properties? Under what conditions is the category of 
$g$-twisted $V$-modules for all
$g\in G$  a $G$-crossed (tensor) category?
\end{prob}

\begin{rema}
{\rm   In the case of $G=\{1\}$, Conjecture \ref{twisted-int-op} is a theorem
(see \cite{H-diff-eqn} and \cite{H-modular}). From this theorem, the author constructed a modular tensor category (see \cite{H-rigidity})
and thus
Conjecture \ref{g-crossed} is also a theorem in this case. }
\end{rema}

In the remaining part of this paper, we describe in details 
the construction and study of twisted modules, the basic and conjectural 
properties of twisted intertwining operators and some thoughts of the author on 
further developments of orbifold conformal field theory. 

\renewcommand{\theequation}{\thesection.\arabic{equation}}
\renewcommand{\thethm}{\thesection.\arabic{thm}}
\setcounter{equation}{0}
\setcounter{thm}{0}
\section{Twisted modules, a general construction and existence results}

To construct orbifold conformal field theories, we first  have to understand  the structures and 
properties of twisted modules. 
 Twisted modules associated to automorphisms of finite orders of vertex operator algebras
 appeared first in the works of  Frenkel-Lepowsky-Meurman \cite{FLM} and Lepowsky \cite{L}. 
In \cite{H-log-twisted-mod}, the author introduced twisted modules associated to 
automorphisms of arbitrary orders (including in particular, infinite orders).
In the case
of automorphisms of infinite orders, 
the logarithm of the variable might appear in twisted vertex operators. Here we 
recall $g$-twisted modules and their variants from \cite{H-log-twisted-mod}.

\begin{defn}
{\rm Let $V$ be a grading-restricted vertex algebra or a vertex operator algebra. 
A generalized  $g$-twisted
$V$-module is a ${\C}$-graded
vector space $W = \coprod_{n \in \C} W_{[n]}$ (graded by {\it weights})
equipped with a linear map
\begin{eqnarray*}
Y_{W}^g: V\otimes W &\to& W\{x\}[\log  x],\\
v \otimes w &\mapsto &Y_{W}^g(v, x)w
\end{eqnarray*}
satisfying the following conditions:

\begin{enumerate}

\item The {\it equivariance property}: For $p \in \mathbb{Z}$, $z
\in \mathbb{C}^{\times}$,  $v \in V$ and $w \in W$, 
$$(Y^{g}_{W})^{p + 1}(gv,
z)w = (Y^{g}_{W})^{p}(v, z)w,$$
where for $p \in \mathbb{Z}$, $(Y^{g}_{W})^{ p}(v, z)$
is the $p$-th analytic branch of $Y_{W}^g(v, x)$.

\item The {\it identity property}: For $w \in W$, $Y^g({\bf 1}, x)w
= w$.

\item  The {\it duality property}: For
any $u, v \in V$, $w \in W$ and $w' \in W'$, there exists a
multivalued analytic function with preferred branch of the form
\[
f(z_1, z_2) = \sum_{i,
j, k, l = 0}^N a_{ijkl}z_1^{m_i}z_2^{n_j}({\rm log}z_1)^k({\rm
log}z_2)^l(z_1 - z_2)^{-t}
\]
for $N \in \mathbb{N}$, $m_1, \dots,
m_N$, $n_1, \dots, n_N \in \mathbb{C}$ and $t \in \mathbb{Z}_{+}$,
such that the series
\begin{align*}
&\langle w', (Y^{g}_{W})^{ p}(u, z_1)(Y^{g}_{W})^{p}(v,
z_2)w\rangle,\\
&\langle w', (Y^{g}_{W})^{ p}(v,
z_2)(Y^{g}_{W})^{p}(u, z_1)w\rangle,\\
&\langle w', (Y^{g}_{W})^{ p}(Y_{V}(u, z_1 - z_2)v,
z_2)w\rangle
\end{align*}
are absolutely convergent in
the regions $|z_1| > |z_2| > 0$, $|z_2| > |z_1| > 0$, $|z_2| > |z_1
- z_2| > 0$, respectively, 
and their sums are equal to the branch
\[
f^{p,p}(z_{1}, z_{2})= \sum_{i, j, k, l = 0}^N
a_{ijkl}e^{m_il_p(z_1)}e^{n_jl_p(z_2)}l_p(z_1)^kl_p(z_2)^l(z_1 -
z_2)^{-t}
\]
of $f(z_1, z_2)$ in the region $|z_1| > |z_2| > 0$, the region $|z_2| > |z_1| > 0$, 
the region given by $|z_2| > |z_1- z_2| > 0$ and $|\arg z_{1}-\arg z_{2}|<\frac{\pi}{2}$, respectively.

\item The {\it $L(0)$-grading condition} and {\it $g$-grading condition}: 
Let  $L_{W}^{g}(0)=\res_{x}xY_{W}^g(\omega, x)$.  Then for $n\in \C$ and $\alpha\in \C/\Z$,
$w \in W_{[n]}^{[\alpha]}$,
there exists $K, \Lambda\in \Z_{+}$ such that $(L_{W}^g (0)-n)^{K} w
=(g-e^{2\pi \alpha i})^{\Lambda}w=0$.

\item The $L(-1)$-{\it derivative property}: For $v \in V$,
\[
\frac{d}{dx}Y^g_{W}(v, x) = Y^g_{W}(L_{V}(-1)v, x).
\]
A {\it lower-bounded generalized $g$-twisted $V$-module} is a generalized 
$g$-twisted $V$-module satisfying the condition $W_{[n]}=0$ for $\Re(n)$ sufficiently 
negative. A {\it grading-restricted generalized $g$-twisted $V$-module} is a 
lower-bounded generalized $g$-twisted $V$-module satisfying in addition the condition 
$\dim W_{[n]}<\infty$ for $n\in \C$. An {\it ordinary $g$-twisted $V$-module} 
or simply a {\it $g$-twisted $V$-module} is a grading-restricted generalized $g$-twisted $V$-module
such that $L_{W}(0)$ acts on $W$ semisimply. 
\end{enumerate}}
\end{defn}

In \cite{Li}, Li studied twisted modules associated to automoprhisms of finite order of a 
vertex operator algebra using the method of weak commutativity and applied this method to 
such twisted modules for vertex operator (super)algebras obtained from infinite-dimensional Lie (super)algebras. 
In \cite{Ba}, Bakalov gave a Jacobi identity for twisted modules associated to 
automorphisms of possibly infinite order, which can be used 
to replace the duality property in the definition above. See \cite{HY} for a proof. 

In the representation theories of associative algebras and of Lie algebras,  free modules
and Verma modules, respectively,  play fundamental roles. For vertex operator algebras, 
though there had been constructions of twisted modules for some examples 
of vertex operator algebras, for more than thirty years, no twisted modules analogous to 
free modules and Verma modules were constructed. In 2019, the author successfully constructed 
 such analogues. These are lower-bounded generalized $g$-twisted $V$-modules 
satisfying a universal property. For the construction in the case that $V$ is a 
grading-restricted vertex algebra, we refer the reader to \cite{H-const-twisted-mod}.
A minor modification of the construction 
in \cite{H-const-twisted-mod} gives such a universal lower-bounded generalized 
twisted $V$-modules in the case that $V$ is a vertex operator algebra (see 
\cite{H-affine-twisted-mod}).
Here we state and discuss the  existence results
and some properties,
including in particular the universal property, of  lower-bounded and grading-restricted 
generalized $g$-twisted $V$-modules
from  \cite{H-const-twisted-mod} and  \cite{H-exist-twisted-mod}.

\begin{thm}[\cite{H-const-twisted-mod}]\label{universal}
Let $M$ be a vector space with actions of $g$ and an operator 
$L_{M}(0)$. Assume that there exist operators $\mathcal{L}_{g}$,
$\mathcal{S}_{g}$, $\mathcal{N}_{g}$ such that on $M$, $g=e^{2\pi i\mathcal{L}_{g}}$
and $\mathcal{S}_{g}$ and $\mathcal{N}_{g}$ are the semisimple and nilpotent, respectively,
parts of $\mathcal{L}_{g}$. Also assume that $L_{M}(0)$ can be decomposed as 
the sum of its semisimple part $L_{M}(0)_{S}$ and nilpotent part $L_{M}(0)_{N}$ and
that the real parts of the eigenvalues of 
$L_{M}(0)$ has a lower bound.   Let $B\in \R$ such that $B$ is less than 
or equal to the real parts of the eigenvalues of $L_{M}(0)$ on $M$. 
Then there exists a lower-bounded generalized $g$-twisted $V$-module 
 $\widehat{M}^{[g]}_{B}$ satisfying the following universal property:
Let $(W, Y^{g}_{W})$ be a lower-bounded generalized $g$-twisted $V$-module
such that $W_{[n]}=0$ when $\Re(n)<B$ and let $M_{0}$ be 
a subspace of $W$ invariant under the actions of 
$g$, $\mathcal{S}_{g}$, $\mathcal{N}_{g}$, $L_{W}(0)$, 
$L_{W}(0)_{S}$ and
$L_{W}(0)_{N}$. 
Assume that there is a linear map $f: M\to M_{0}$  commuting with the actions of 
$g$, $\mathcal{S}_{g}$, $\mathcal{N}_{g}$, $L_{M}(0)$ and $L_{W}(0)|_{M_{0}}$, 
$L_{M}(0)_{S}$ and $L_{W}(0)_{S}|_{M_{0}}$ and
$L_{M}(0)_{N}$ and $L_{W}(0)_{N})|_{M_{0}}$. Then there exists a unique module 
map $\hat{f}: \widehat{M}^{[g]}_{B}\to W$ such that $\hat{f}|_{M}=f$. 
If $f$ is surjective and $(W, Y^{g}_{W})$ is generated by the 
coefficients of $(Y^{g})_{WV}^{W}(w_{0}, x)v$ for $w_{0}\in M_{0}$ and $v\in V$, 
where $(Y^{g})_{WV}^{W}$ 
is the twist vertex operator map obtained from $Y_{W}^{g}$, then 
$\hat{f}$ is surjective. 
\end{thm}

An explicit construction, not just the existence, of $\widehat{M}^{[g]}_{B}$
was given in \cite{H-const-twisted-mod} for a grading-restricted vertex algebra $V$. 
One crucial ingredient in this construction is the twist vertex operators introduced and 
studied in \cite{H-twist-vo}.
In the case that $V$ is a vertex operator algebra, see Subsection 4.1 of \cite{H-affine-twisted-mod}.
The lower-bounded generalized $g$-twisted $V$-module 
 $\widehat{M}^{[g]}_{B}$ is also unique up to equivalence by the universal property. 

One immediate consequence of Theorem \ref{universal} is the following result:

\begin{cor}[\cite{H-const-twisted-mod}]
Let $(W, Y^{g}_{W})$ be a lower-bounded generalized $g$-twisted $V$-module generated by the 
coefficients of $(Y^{g})_{WV}^{W}(w, x)v$ for $w\in M$, where $(Y^{g})_{WV}^{W}$ 
is the twist vertex operator map obtained from $Y_{W}^{g}$
and $M$ is a $\Z_{2}$-graded subspace of $W$  invariant under the actions of 
$g$, $\mathcal{S}_{g}$, $\mathcal{N}_{g}$, $L_{W}(0)$, 
$L_{W}(0)_{S}$ and
$L_{W}(0)_{N}$. Let $B\in \R$ such that $W_{[n]}=0$ when $\Re(n)<B$. 
Then there is a generalized $g$-twisted $V$-submodule $J$ of 
$\widehat{M}^{[g]}_{B}$ such that $W$ is equivalent as a 
lower-bounded generalized $g$-twisted $V$-module to the quotient 
module $\widehat{M}^{[g]}_{B}/J$.
\end{cor}

The construction of $\widehat{M}^{[g]}_{B}$ has been used 
in \cite{H-exist-twisted-mod} to solve some
open problems of more than twenty years. We now discuss these results.
These results are proved in \cite{H-exist-twisted-mod} when $V$ is a 
grading-restricted vertex algebra or a M\"{o}bius vertex algebra. But 
by using the minor modification in \cite{H-affine-twisted-mod} of the construction of 
universal lower-bounded generalized $g$-twisted $V$-modules
 in the case that $V$ is a vertex operator algebra, the same proofs also work in this case. 

The first  such open problem is the existence of nonzero lower-bounded generalized 
$g$-twisted $V$-modules. Note that since $V$ itself  is a $V$-module, the existence of 
nonzero $V$-modules is obvious. But it is highly nontrivial why nonzero lower-bounded 
generalized $g$-twisted $V$-modules exist. Assuming that the vertex operator 
algebra is simple and $C_{2}$-cofinite and the automorphism
is of finite order, 
Dong, Li and Mason proved the existence of an
irreducible twisted module \cite{DLM2}. But no progress has been made in the 
general case for more than twenty years. 
On the other hand, some classes of vertex operator algebras
that are not $C_{2}$-cofinite have a very rich and 
exciting representation theory. For example,  vertex operator algebras
associated to affine Lie algebras at admissible levels are not $C_{2}$-cofinite. 
But the category of ordinary modules 
for such a vertex operator algebra has a braided tensor category structure with a twist (see \cite{CHY}),
which is also rigid  and
in many cases even has a modular tensor category 
structure (see \cite{CHY} for the conjectures and a proof in the case of $\mathfrak{sl}_{2}$ and  see \cite{C}
for a proof of the rigidity in the simply-laced case).  It is important to study 
vertex operator algebras that are not $C_{2}$-cofinite. 
 From Theorem \ref{universal}, fora  general grading-restricted vertex algebra $V$
(not necessarily $C_{2}$-cofinite)
we indeed have constructed lower-bounded generalized 
$g$-twisted $V$-modules. But 
it is not obvious from the construction in 
\cite{H-const-twisted-mod} why $\widehat{M}^{[g]}_{B}$
is not $0$. In \cite{H-exist-twisted-mod}, the author solved this problem completely.

\begin{thm}[\cite{H-exist-twisted-mod}]\label{existence}
The lower-bounded generalized $g$-twisted $V$-module 
$\widehat{M}^{[g]}_{B}$ is not $0$. In particular, there exists 
nonzero lower-bounded generalized $g$-twisted $V$-modules. 
\end{thm}

In \cite{DLM1}, Dong, Li and Mason generalized Zhu's algebra $A(V)$ (see \cite{Z}) to 
a twisted Zhu's algebra $A_{g}(V)$ for a vertex operator algebra $V$
and an automorphism $g$ of $V$ of finite order.
In \cite{HY}, Yang and the author introduced twisted zero-mode algebra $Z_{g}(V)$ associated to 
$V$ and an automorphism $g$ of $V$ not necessarily of finite order and also generalized 
the twisted Zhu's algebra $A_{g}(V)$ to the case that $g$ is not of finite order. These two associative algebras
are in fact isomorphic (see \cite{HY}). In the case that $g$ is of finite order, 
$A_{g}\ne 0$ is stated explicitly as a conjecture  in the beginning of Section 9 of the arXiv version 
of \cite{DLM1}. 
In the case that $V$ is $C_{2}$-cofinite and $g$ is of finite order, Dong, Li and Mason proved
this conjecture  in \cite{DLM2}. But in general, this conjecture had been open until 
it was proved in \cite{H-exist-twisted-mod} as an immediate consequence of Theorem \ref{existence}. 

\begin{cor}[\cite{H-exist-twisted-mod}]
The twisted Zhu's algebra $A_{g}(V)$ or the twisted zero-mode algebra 
$Z_{g}(V)$ is not $0$. 
\end{cor}

Another application of Theorem \ref{existence} is on the following existence of irreducible 
lower-bounded  generalized $g$-twisted $V$-module:

\begin{thm}[\cite{H-exist-twisted-mod}]
Let $W$ be a lower-bounded generalized $g$-twisted $V$-module generated by 
a nonzero element $w$ (for example, $\widehat{M}^{[g]}_{B}$
when $M$ is a one dimensional space spanned by an element $w$ and 
$B$ is less than or equal to the real part of the weight 
of $w$). 
Then there exists a maximal submodule $J$ of $W$ such that 
$J$ does not contain $w$ and the quotient 
$W/J$ is irreducible.
\end{thm}

Though lower-bounded generalized $g$-twisted $V$-modules are important in our study of orbifold conformal
field theories, we are mainly interested in grading-restricted generalized $g$-twisted $V$-modules
and ordinary $g$-twisted $V$-modules,
because their double contragredients  are equivalent to themselves. 
One  important problem is the existence of irreducible grading-restricted generalized $g$-twisted $V$-modules
and irreducible ordinary $g$-twisted $V$-modules.
In the case that $V$ is simple, $C_{2}$-cofinite and $g$ is of finite order, Dong, Li and Mason 
proved in \cite{DLM2}  the existence of an irreducible 
ordinary $g$-twisted  $V$-module. Their proof used genus-one $1$-point functions. Thus  the simplicity and 
$C_{2}$-cofiniteness of $V$ and the finiteness of the order of $g$ are necessary in their approach. 
Using our construction of the universal lower-bounded generalized $g$-twisted $V$-modules,
the author proved in \cite{H-exist-twisted-mod} the existence of  irreducible grading-restricted generalized
$g$-twisted  $V$-modules and irreducible ordinary $g$-twisted  $V$-modules
under some very weak conditions. In particular, 
the simplicity and $C_2$-cofininess of $V$ and the finiteness of the order of $g$ are not needed. 

\begin{thm}[\cite{H-exist-twisted-mod}]
Let $V$ be a M\"{o}bius vertex superalgebra
and $g$ an automorphism of $V$.  Assume that the set
of real parts of the numbers in $P(V)$ has no accumulation point in $\R$. 
If the twisted Zhu's algebra $A_{g}(V)$ or the twisted zero-mode algebra $Z_{g}(V)$
is finite dimensional, then there exists an irreducible grading-restricted generalized
$g$-twisted  $V$-module. Such an irreducible grading-restricted generalized
$g$-twisted  $V$-module is an irreducible ordinary $g$-twisted  $V$-module if $g$ acts 
on it semisimply. In particular, if $g$ is of finite order, there exists an irreducible 
ordinary $g$-twisted  $V$-module.
\end{thm}

The author also proved in \cite{H-exist-twisted-mod} that a lower-bounded generalized 
module for the fixed-point subalgebra $V^{g}$ of $V$ can be extended to a 
 lower-bounded generalized $g$-twisted $V$-module. 

\begin{thm}[\cite{H-exist-twisted-mod}]\label{v-g-extension}
Let $V$ be a grading-restricted vertex algebra and 
$W_{0}$ a lower-bounded generalized $V^{g}$-module (in particular, 
$W_{0}$ has a lower-bounded grading by $\C$).
Assume that 
$g$ acts on $W_{0}$ and there are semisimple and nilpotent operators $\mathcal{S}_{g}$
and $\mathcal{N}_{g}$, respectively, on $W_{0}$ such that 
$g=e^{2\pi i\mathcal{L}_{g}}$ where $\mathcal{L}_{g}=\mathcal{S}_{g}+\mathcal{N}_{g}$. 
Then $W_{0}$ can be extended to a lower-bounded generalized $g$-twisted $V$-module,
that is,  there exists a lower-bounded generalized $g$-twisted $V$-module $W$ 
and an injective module map $f: W_{0}\to W$ of $V^{g}$-modules. 
\end{thm}

We have the following open problem:

\begin{prob}
 Finding conditions on the vertex operator algebra such that under these conditions, 
 the universal lower-bounded generalized $g$-twisted $V$-modules
have irreducible quotients whose homogeneous subspaces
are finite dimensional. 
\end{prob}

\begin{rema}
{\rm In the case that $V$ is a grading-restricted vertex algebra or vertex operator algebra
associated to an affine Lie algebra, the author solved this problem (see \cite{H-affine-twisted-mod}). 
Since the results in \cite{H-affine-twisted-mod} are only for special types of examples of grading
restricted vertex algebras or vertex operator algebras, we shall not discuss 
the details here. The interested reader is referred to the paper \cite{H-affine-twisted-mod} for details.}
\end{rema}

\renewcommand{\theequation}{\thesection.\arabic{equation}}
\renewcommand{\thethm}{\thesection.\arabic{thm}}
\setcounter{equation}{0}
\setcounter{thm}{0}
\section{Twisted intertwining operators}

To construct orbifold conformal field theories, one approach is to study the 
representation theory of fixed point vertex operator algebras. 
If one can prove that intertwining operators among suitable  modules for the fixed-point vertex operator 
algebra of a vertex operator algebra under a group of automorphisms
satisfy the convergence and extension property, the associativity (operator product expansion),
the modular invariance property and the higher-genus convergence property as is decribed in Section 2, 
then one can construct and study the corresponding orbifold conformal field theories 
using the steps described in Section 2. But it is itself a difficult problem to prove these properties 
using  the properties of 
intertwining operators among modules for the larger vertex operator algebra. 

On the other hand, as is mentioned in Section 3, since twisted modules for the larger vertex operator algebra 
are analogues of (untwisted) modules for the vertex operator algebra, we expect that the properties of 
intertwining operators among twisted modules can be studied by generalizing the results and 
approach for the intertwining operators among modules. Moreover, we expect 
that every module for the fixed-point vertex operator algebra can be obtained 
from a twisted module for the larger vertex operator algebra. In particular, 
intertwining operators among modules for the fixed-point vertex operator algebra can also be obtained from intertwining 
operators among twisted modules for the larger vertex operator algebra. Therefore 
instead of studying intertwining operators among modules for the fixed-point vertex operator algebra, 
we study intertwining operators among twisted modules for the larger vertex operator algebra. 
For simplicity, as in \cite{H-twisted-int},
we call intertwining operators among twisted modules twisted intertwining operators. 

We still need a precise definition of twisted intertwining operator. 
By generalizing the Jacobi identity for intertwining operators, 
Xu introduced the notion of intertwining operators among twisted modules associated 
to commuting automorphisms of finite orders (see \cite{X}). 
But in general, an orbifold conformal field theory might be associated to a nonabelian group of automorphisms.
In particular, we have to introduce and study intertwining operators among twisted modules associated 
to not-necessarily-commuting automorphisms. Also, the group might not be finite. So we also have to introduce and study 
intertwining operators among twisted modules associated 
to automorphisms of infinite orders. 
The formulation used in \cite{FHL} and \cite{X} cannot be generalized directly to give a definition of 
intertwining operators among twisted modules associated 
to not-necessarily-commuting automorphisms. 
For more than twenty years, no such definition was given in the literature. 
This is the reason why orbifold conformal field theories associated to nonabelian groups had not been studied much
mathematically in the past. 

 In \cite{H-twisted-int}, the author found a formulation of such a notion of 
twisted intertwining operators associated 
to not-necessarily-commuting automorphisms of possibly infinite orders and proved their basic properties. 
The general theory and construction of orbifold conformal field theories associated to
nonabelian groups (including infinite groups) can now be started from such operators. 

We first give the precise definition of twisted intertwining operators.

\begin{defn}
Let $g_{1}, g_{2}, g_{3}$ be automorphisms of $V$ and 
let $W_{1}$, $W_{2}$ and $W_{3}$ be $g_{1}$-, $g_{2}$- 
and $g_{3}$-twisted 
$V$-modules, respectively. A {\it twisted intertwining operator of type $\binom{W_{3}}{W_{1}W_{2}}$} is a 
linear map
\begin{eqnarray*}
\mathcal{Y}: W_{1}\otimes W_{2}&\to& W_{3}\{x\}[\log x]\nn
w_{1}\otimes w_{2}&\mapsto& \mathcal{Y}(w_{1}, x)w_{2}=\sum_{k=0}^{K}\sum_{n\in \C}
\mathcal{Y}_{n, k}(w_{1})w_{2}x^{-n-1}(\log x)^{k}
\end{eqnarray*}
satisfying the following conditions:

\begin{enumerate}

\item {\it The lower truncation property}: For $w_{1}\in W_{1}$ and $w_{2}\in W_{2}$, $n\in \C$ and $k=0, \dots, K$, 
$\mathcal{Y}_{n+l, k}(w_{1})w_{2}=0$ for $l\in \N$ sufficiently large.

\item The {\it duality property}:  For $u\in V$, $w_{1}\in W_{1}$, $w_{2}\in W_{2}$
and $w_{3}'\in W_{3}'$, there exists a
multivalued analytic function with preferred branch
\begin{align*}
f&(z_1, z_2; u, w_{1}, w_{2}, w_{3}')\nn
&\quad 
= \sum_{i, j, k, l, m, n=1}^{N}a_{ijklmn}z_1^{r_{i}}z_2^{s_{j}}(z_1 - z_2)^{t_{k}}(\log z_1)^l
(\log z_2)^{m}(\log (z_1 - z_2))^{n}
\end{align*}
for $N \in \mathbb{N}$, $r_{i}, s_{j}, t_{k}, a_{ijklmn}\in \mathbb{C}$,
such that for $p_{1}, p_{2}, p_{12}\in \Z$, the series
\begin{align*}
&\langle w'_{3}, (Y_{W_{3}}
^{g_{3}})^{p_{1}}(u, z_{1})\mathcal{Y}^{p_{2}}(w_{1}, z_{2})w_{2}\rangle,\\
&\langle w'_{3}, \mathcal{Y}^{p_{2}}(w_{1}, z_{2})
(Y_{W_{2}}^{g_{2}})^{p_{1}}(u, z_{1})w_{2}\rangle,\\
& \langle w'_{3}, \mathcal{Y}^{p_{2}}(
(Y_{W_{1}}^{g_{1}})^{p_{12}}(u, z_{1}-z_{2})v, z_{2})w\rangle
\end{align*}
are absolutely convergent in the regions 
$|z_{1}|>|z_{2}|>0$,
$|z_{2}|>|z_{1}|>0$, $|z_{2}|>|z_{1}-z_{2}|>0$, respectively. Moreover, their sums are equal 
to the branches
\begin{align*}
f&^{p_{1}, p_{2}, p_{1}}(z_{1}, z_{2}; u, w_{1}, w_{2}, w_{3}')\nn
&=\sum_{i, j, k, l, m, n=1}^{N}a_{ijklmn}z_1^{r_{i}}e^{r_{i}l_{p_{1}}(z_{1})}
e^{s_{j}l_{p_{2}}(z_2)}e^{t_{k}l_{p_{1}}(z_1 - z_2)} (l_{p_{1}}(z_1))^l
(l_{p_{2}}(z_2))^{m}(l_{p_{1}}(z_1 - z_2))^{n},\nn
f&^{p_{1}, p_{2}, p_{2}}(z_{1}, z_{2}; u, w_{1}, w_{2}, w_{3}')\nn
&=\sum_{i, j, k, l, m, n=1}^{N}a_{ijklmn}z_1^{r_{i}}e^{r_{i}l_{p_{1}}(z_{1})}
e^{s_{j}l_{p_{2}}(z_2)}e^{t_{k}l_{p_{2}}(z_1 - z_2)}(l_{p_{1}}(z_1))^l
(l_{p_{2}}(z_2))^{m}(l_{p_{2}}(z_1 - z_2))^{n},\nn
f&^{p_{2}, p_{2}, p_{12}}(z_{1}, z_{2}; u, w_{1}, w_{2}, w_{3}')\nn
&=\sum_{i, j, k, l, m, n=1}^{N}a_{ijklmn}z_1^{r_{i}}e^{r_{i}l_{p_{2}}(z_{1})}
e^{s_{j}l_{p_{2}}(z_2)}e^{t_{k}l_{p_{12}}(z_1 - z_2)} (l_{p_{2}}(z_1))^l
(l_{p_{2}}(z_2))^{m}(l_{p_{12}}(z_1 - z_2))^{n},
\end{align*}
respectively, of $f(z_1, z_2; u, w_{1}, w_{2}, w_{3}')$ in the region given by 
$|z_{1}|>|z_{2}|>0$ and $|\arg (z_{1}-z_{2})-\arg z_{1}|<\frac{\pi}{2}$,
the region given by $|z_{2}|>|z_{1}|>0$ and $-\frac{3\pi}{2}< \arg (z_{1}-z_{2})-\arg z_{2}<-\frac{\pi}{2}$,
the region given by $|z_{2}|>|z_{1}-z_{2}|>0$ and $|\arg  z_{1}- \arg z_{2}|<\frac{\pi}{2}$, respectively.

\item The {\it $L(-1)$-derivative property}: 
$$\frac{d}{dx}\mathcal{Y}(w_{1}, x)=\mathcal{Y}(L(-1)w_{1}, x).$$

\end{enumerate}
\end{defn}

The correct notion of twisted intertwining operator should have some 
basic properties. These properties were proved in \cite{H-twisted-int}.
In particular, the notion of twisted intertwining operator introduced in 
\cite{H-twisted-int} is  indeed the correct one. 

The first property is the following: 

\begin{thm}[\cite{H-twisted-int}]
 Let $g_{1}, g_{2}, g_{3}$ be automorphisms of $V$ and let $W_{1}$, $W_{2}$ and $W_{3}$ be $g_{1}$-, $g_{2}$- 
and $g_{3}$-twisted 
$V$-modules, respectively.
Assume that the vertex operator map for $W_{3}$ given by $u\mapsto Y_{W_{3}}^{g_{3}}(u, x)$
is injective. If there exists a twisted intertwining operator $\mathcal{Y}$ of type $\binom{W_{3}}{W_{1}W_{2}}$ such that 
the coefficients of the series $\mathcal{Y}(w_{1}, x)w_{2}$ for $w_{1}\in W_{1}$ and $w_{2}\in W_{2}$ span
$W_{3}$,
then $g_{3}=g_{1}g_{2}$.
\end{thm}

For the proof of this theorem, we refer the reader to \cite{H-twisted-int}.
Here we give a geometric explanation of the theorem  in terms of a picture (Figure 1). 
From the equivariance property for twisted modules, the monodromy of 
the twisted vertex operators corresponds to the action of the automorphism
on $V$. Given a twisted intertwining operator of type 
$\binom{W_{3}}{W_{1}W_{2}}$, the monodromy of the twisted vertex operator
for $W_{3}$ gives $g_{3}$. This is described in the left braiding graph in 
Figure 1. But this braiding graph is topologically the same as the right braiding 
graph in Figure 1. It is clear that the right braiding graph in Figure 1 is equal to the 
product of the monodromy of the twisted vertex operator for $W_{1}$ and 
the monodromy of the twisted vertex operator for $W_{2}$. So the 
right braiding graph in Figure 1 gives $g_{1}g_{2}$. Thus we see from this 
geometric picture, $g_{3}$ should be equal to $g_{1}g_{2}$.

\begin{figure}
\begin{center}
\includegraphics
{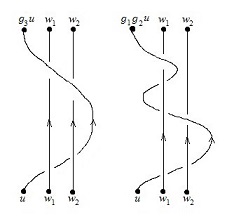}
\end{center}
\caption{The braiding graphs corresponding to $g_{3}$ (left) and 
$g_{1}g_{2}$ (right)\label{figure3}}
\end{figure}

For the other properties, we first need to recall an action of an automorphism 
$h$ of $V$ on a $g$-twisted $V$-module. Let $(W, Y^{g}_{W})$ be a  $g$-twisted
$V$-module. Let $h$ be an automorphism of $V$ and let
\begin{eqnarray*}
\phi_{h}(Y^{g}): V\times W&\to& W\{x\}[{\rm log} x]\nn
v \otimes w & \mapsto & \phi_{h}(Y^{g})(v, x)w
\end{eqnarray*}
be the linear map defined by
$$ \phi_{h}(Y^{g})(v, x)w=Y^{g}(h^{-1}v, x)w.$$
Then the pair $(W, \phi_{h}(Y^{g}))$ is an $hgh^{-1}$-twisted
$V$-module. We shall denote the $hgh^{-1}$-twisted
$V$-module in the proposition above by $\phi_{h}(W)$.

We now discuss the skew-symmetry isomorphism for twisted intertwining operators. 
Let $g_{1}, g_{2}$ be automorphisms of $V$,  $W_{1}$, $W_{2}$ and $W_{3}$ $g_{1}$-, $g_{2}$- 
and $g_{1}g_{2}$-twisted 
$V$-modules and $\mathcal{Y}$ a twisted intertwining operator 
of type $\binom{W_{3}}{W_{1}W_{2}}$. 
We define  linear maps
\begin{eqnarray*}
\Omega_{\pm}(\Y): W_{2}\otimes W_{1}&\to& W_{3}\{x\}[\log x]\nn
w_{2}\otimes w_{1}&\mapsto& \Omega_{\pm}(\Y)(w_{2}, x)w_{1}
\end{eqnarray*}
by
$$
\Omega_{\pm}(\Y)(w_{2}, x)w_{1}=e^{xL(-1)}\Y(w_{1}, y)w_{2}\lbar_{y^{n}=e^{\pm n\pi i}x^{n}, \;
\log y=\log x\pm \pi i}
$$
for $w_{1}\in W_{1}$ and $w_{2}\in W_{2}$. 

\begin{thm}[\cite{H-twisted-int}]\label{omega}
The linear maps $\Omega_{+}(\Y)$ and $\Omega_{-}(\Y)$  are twisted intertwining operators
of types $\binom{W_{3}}{W_{2}\;\phi_{g_{2}^{-1}}(W_{1})}$ and $\binom{W_{3}}{\phi_{g_{1}}(W_{2})\;W_{1}}$,
respectively. 
\end{thm}

From Theorem \ref{omega}, we see that $\Omega_{+}$ and $\Omega_{-}$ are indeed 
ismorphisms between spaces of twisted intertwining operators. 

\begin{cor}[\cite{H-twisted-int}]
The maps $\Omega_{+}: \mathcal{V}_{W_{1}W_{2}}^{W_{3}}\to \mathcal{V}_{W_{2}\phi_{g_{2}^{-1}}(W_{1})}^{W_{3}}$ and 
$\Omega_{-}: \mathcal{V}_{W_{1}W_{2}}^{W_{3}}\to \mathcal{V}_{\phi_{g_{1}}(W_{2})W_{1}}^{W_{3}}$
are linear isomorphisms. In particular, 
$\mathcal{V}_{W_{1}W_{2}}^{W_{3}}$,  $\mathcal{V}_{\phi_{g_{1}}(W_{2})W_{1}}^{W_{3}}$
and $\mathcal{V}_{W_{2}\phi_{g_{2}^{-1}}(W_{1})}^{W_{3}}$
are linearly isomorphic.
\end{cor}

Finally we discuss the contragredient isomorphism for twisted intertwining operators. 
We first recall contragredient twisted $V$-modules. 
Let $(W, Y^{g}_{W})$ be a $g$-twisted $V$-module relative to $G$.
Let $W'$ be the graded dual of $W$. Define
a linear map 
\begin{eqnarray*}
(Y_{W}^{g})': V\otimes W' &\to& W'\{x\}[{\rm log} x],\\
v \otimes w' &\mapsto &(Y^g_{W})'(v, x)w'
\end{eqnarray*}
by
$$\langle (Y^g_{W})'(v, x)w', w\rangle=\langle w', Y^{g}_{W}(e^{xL(1)}(-x^{-2})^{L(0)}v, x^{-1})w\rangle$$
for $v\in V$, $w\in W$ and $w'\in W'$. 
Then the pair $(W', (Y^{g}_{W})')$ is a $g^{-1}$-twisted $V$-module.
We call $(W', (Y^{g}_{W})')$ the contragredient twisted $V$-module of $(W, Y^{g}_{W})$.

Let $g_{1}, g_{2}$ be automorphisms of $V$,  $W_{1}$, $W_{2}$ and $W_{3}$ $g_{1}$-, $g_{2}$- 
and $g_{1}g_{2}$-twisted 
$V$-modules and $\mathcal{Y}$ a twisted intertwining operator 
of type $\binom{W_{3}}{W_{1}W_{2}}$. 
We define  linear maps
\begin{eqnarray*}
A_{\pm}(\Y): W_{1}\otimes W_{3}'&\to& W_{2}'\{x\}[\log x]\nn
w_{1}\otimes w_{3}'&\mapsto& A_{\pm}(\Y)(w_{1}, x)w_{3}'
\end{eqnarray*}
by
\begin{align*}
\langle A_{\pm}(\Y)(w_{1}, x)w_{3}', w_{2}\rangle
&=\langle w_{3}', \Y(e^{xL(1)}e^{\pm \pi i L(0)}(x^{-L(0)})^{2}w_{1}, x^{-1})w_{2}\rangle
\end{align*}
for $w_{1}\in W_{1}$ and $w_{2}\in W_{2}$ and $w_{3}'\in W_{3}'$.

\begin{thm}[\cite{H-twisted-int}]\label{contragredient}
The linear maps $A_{+}(\Y)$ and $A_{-}(\Y)$  are twisted intertwining operators
of types $\binom{\phi_{g_{1}}(W_{2}')}{W_{1}W_{3}'}$ and $\binom{W_{2}'}{W_{1}\phi_{g_{1}^{-1}}(W_{3}')}$, respectively.
\end{thm}

From Theorem \ref{contragredient}, we see that $A_{+}$ and $A_{-}$ are indeed 
ismorphisms between spaces of twisted intertwining operators. 

\begin{cor}[\cite{H-twisted-int}]
The maps $A_{+}: \mathcal{V}_{W_{1}W_{2}}^{W_{3}}\to \mathcal{V}_{W_{1}W_{3}'}^{\phi_{g_{1}}(W_{2}')}$ and 
$A_{-}: \mathcal{V}_{W_{1}W_{2}}^{W_{3}}\to \mathcal{V}_{W_{1}\phi_{g_{1}^{-1}}(W_{3}')}^{W_{2}'}$
are linear isomorphisms. In particular, 
$\mathcal{V}_{W_{1}W_{2}}^{W_{3}}$,  $\mathcal{V}_{W_{1}W_{3}'}^{\phi_{g_{1}}(W_{2}')}$
and $\mathcal{V}_{W_{1}\phi_{g_{1}^{-1}}(W_{3}')}^{W_{2}'}$
are linearly isomorphic.
\end{cor}

\renewcommand{\theequation}{\thesection.\arabic{equation}}
\renewcommand{\thethm}{\thesection.\arabic{thm}}
\setcounter{equation}{0}
\setcounter{thm}{0}
\section{Main conjectural properties of twisted intertwining operators}

In Section 3, we have recalled some main conjectures on the construction of 
orbifold conformal field theories using the representation theory of vertex operator algebras. 
In this section, we formulate precisely and discuss  in details the conjectural properties in these conjectures.

In several conjectures and problems in Section 3, only associativity, commutativity and modular invariance of
twisted intertwining operators are stated. But as in the untwisted case, the statements of these properties 
make sense only when the corresponding convergence properties hold. In fact, the proofs of these 
convergence properties are one of the main difficult parts of the proofs of these properties in the 
untwisted case. The twisted case will certainly be the same. So below we shall first give the formulation of these 
convergence together with the analytic extensions of the convergent series before the formulations 
of these properties themselves. 

We first formulate the conjectural convergence and extension property, 
associativity and commutativity of twisted intertwining operators. 
We formulate them for lower-bounded generalized twisted $V$-modules so that 
they are more flexible. 
But we warn the reader that these properties in general will not be true for such general twisted modules. 
Usually the twisted modules should be grading-restricted and satisfy some additional conditions. 
In the case that the vertex operator algebra is finite reductive and 
 the group of automorphisms is finite, the category of twisted modules for which 
these properties holds is conjectured to be the category of ordinary twisted modules. 
In general the correct categories of twisted modules will be given in precise conjectures in the future
and are also an important part of the research on the construction of orbifold conformal field theories. 

\vspace{1em}

\noindent {\bf Convergence and extension property of products of $n$ twisted intertwining operators}: 
Let  $g_{1}, \dots, g_{n+1}$
be automorphisms of $V$. Let 
$W_{0}, W_{1},\dots W_{n+1}$,  $\widetilde{W}_{1}, \dots, \widetilde{W}_{n-1}$ be
lower-bounded generalized  $(g_{1}\cdots g_{n+1})$-, $g_{1}$-, 
$\dots$, $g_{n+1}$-, $(g_{2}\cdots g_{n+1})$-, $\dots$, $(g_{n}g_{n+1})$-twisted $V$-modules, 
respectively,  and  $\Y_{1},\dots, \Y_{i}, \dots, \Y_{n}$ twisted 
intertwining operators of types 
$\binom{W_{0}}{W_{1}\;\widetilde{W}_{1}}$, $\dots$, $\binom{\widetilde{W}_{i-1}}{W_{i}\;\widetilde{W}_{i}}$, 
$\dots$, $\binom{\widetilde{W}_{n-1}}{W_n\;W_{n+1}}$, respectively.
For $w_{1}\in W_{1},\ldots, w_{n+1}\in W_{n+1}$  and $w_{0}'\in W_{0}'$, 
the series 
$$\langle w_{0}',\Y_{1}(w_{1},z_{1})\cdots \Y_{n}(w_{n},z_{n})w_{n+1}\rangle$$
in complex variables $z_1,\ldots z_n$  is absolutely convergent in the 
region $|z_1|>\cdots>|z_n|>0$ and its sum can be analytically 
continued to a multivalued analytic function 
$$F(\langle u_{1}',\Y_1(w_1,z_1)\cdots \Y_n(w_n,z_n)u_{n+1}\rangle)$$
on the region
 $$\{(z_1,\ldots,z_n)\mid z_i\ne 0, z_i-z_j\ne 0 \;{\rm for}\; i\ne j\}\subset \C^n$$
 and the only possible singular points $z_i=0,\infty$ and $z_i=z_j$ are regular singular points. 

\vspace{1em}

To formulate the associativity of twisted intertwining operators, we need the following 
result:

\begin{prop}
Assume that the convergence and extension property of products of $2$ twisted intertwining operators holds. 
Let $g_{1}, \dots, g_{4}, g$ be automorphisms of $V$. Let
$W_1,  W_{2}, W_{3}$, $W_{4}, W$ be lower-bounded generalized 
$g_{1}$-, $g_{2}$-, $g_{3}$-, $(g_{1}g_{2}g_{3})$-, $(g_{1}g_{2})$-twisted  $V$-modules
and $\Y_{3}$ and $\Y_{4}$ twisted intertwining operators  of types 
$\binom{W_{4}}{W_{1}\; W}$ and $\binom{W}{W_{2} \; W_{3}}$, respectively.
Then for $w_1\in W_1, w_{2}\in W_{2}$, $w_{3}\in W_{3}$ and $w_{4}'\in W_{4}'$,
the series
$$\langle w_{4}', \Y_3(\Y_{4}(w_1, z_1-z_{2})w_2, z_2)w_{3}\rangle$$
is absolutely convergent in the region $|z_{2}|>|z_{1}-z_{2}|>0$ 
and its sum can be analytically 
continued to a multivalued analytic function 
$$F(\langle w_{4}', \Y_3(\Y_{4}(w_1, z_1-z_{2})w_2, z_2)w_{3}\rangle)$$
on the region
 $$\{(z_1, z_{2})\mid z_1, z_{2}, z_1-z_2\ne 0\}\subset \C^2$$
 and the possible singular points $z_1, z_{2}, z_1-z_2=0, \infty$ are regular singular points. 
\end{prop}

The proof of this proposition is completely the same as the proofs of 
 Proposition 14.1 in \cite{tensor4} and Proposition 7.3 in \cite{HLZ5}. 

We are ready to state precisely the conjectural associativity or the
operator product expansion
of twisted intertwining operators.

\vspace{1em}
\noindent {\bf Associativity  of twisted intertwining operators}
Let  $g_{1},g_{2}, g_{3}$ be automorphisms of $V$. Let
$W_{1}$, $W_2$, $ W_{3}$, $W_{4}$, $W_{5}$ be lower-bounded generalized 
$g_{1}$-, $g_{2}$-, $g_{3}$-, $(g_{1}g_{2}g_{3})$, $(g_{2}g_{3})$-twisted $V$-modules
and $\Y_{1}$ and $\Y_{2}$ intertwining operators of types 
$\binom{W_{4}}{ W_{1}\; W_{5}}$ and $\binom{W_{5}}{W_{2}\; W_{3}}$, respectively. 
There exist a lower-bounded $(g_{1}g_{2})$-twisted generalized $V$-module $W_{6}$
and intertwining operators $\Y_{3}$ and $\Y_{4}$ of the types
$\binom{W_{4}}{ W_{6} \; W_{3}}$ and $\binom{W_{6}}{ W_{1}\;  W_{2}}$, respectively,
such that for $w_1\in W_1, w_{2}\in W_{2}$, $w_{3}\in W_{3}$ and $w_{4}'\in W_{4}'$,
$$F(\langle w_{4}', \Y_1(w_1,z_1) \Y_2(w_2,z_2)w_{3}\rangle)
=F(\langle w_{4}', \Y_3(\Y_{4}(w_1, z_1-z_{2})w_2, z_2)w_{3}\rangle).$$

\vspace{1em}

Another important conjectural property following immediately from the associativity 
of twisted intertwining operators and Theorem \ref{omega}
 is the commutativity of twisted intertwining 
operators:

\vspace{1em}
\noindent {\bf Commutativity  of twisted intertwining operators}
Let  $g_{1},g_{2}, g_{3}$ be automorphisms of $V$. Let
$W_{1}$, $W_2$, $ W_{3}$, $W_{4}$, $W_{5}$ be lower-bounded generalized 
$g_{1}$-, $g_{2}$-, $g_{3}$-, $(g_{1}g_{2}g_{3})$, $(g_{2}g_{3})$-twisted $V$-modules
and $\Y_{1}$ and $\Y_{2}$ intertwining operators of types 
$\binom{W_{4}}{ W_{1} \; W_{5}}$ and $\binom{W_{5}}{W_{2}\; W_{3}}$, respectively. 
There exist a lower-bounded $g_{1}g_{3}$-twisted generalized $V$-module $W_{6}$
and intertwining operators $\Y_{3}$ and $\Y_{4}$ of the types
$\binom{W_{4}}{ \phi_{g_{1}}(W_{2}) \; W_{6}}$ and $\binom{W_{6}}{W_{1}\; W_{3}}$, respectively,
or a  lower-bounded $(g_{2}^{-1}g_{1}g_{2}g_{3})$-twisted generalized $V$-module $W_{6}$
and intertwining operators $\Y_{3}$ and $\Y_{4}$ of the types
$\binom{W_{4}}{ W_{2} \; W_{6}}$ and $\binom{W_{6}}{ \phi_{g_{2}^{-1}}(W_{1})\;  W_{3}}$, respectively,
such that for $w_1\in W_1, w_{2}\in W_{2}$, $w_{3}\in W_{3}$ and $w_{4}'\in W_{4}'$,
$$F(\langle w_{4}', \Y_1(w_1,z_1) \Y_2(w_2,z_2)w_{3}\rangle)
=F(\langle w_{4}', \Y_3(w_2, z_2)\Y_{4}(w_1, z_{1})w_{3}\rangle).$$

\vspace{1em}

The convergence and extension property, the associativity and commutativity 
of twisted intertwining operators discussed above are genus-zero properties.
We now discuss the conjectural genus-one properties. For a conformal field theory, genus-one 
correlation functions should be equal to the analytic extensions of $q$-traces or more generally 
pseudo-$q$-traces of products of geometrically-modified intertwining operators. So in our case, 
we have to consider $q$-trace or pseudo-$q$-traces of products of geometrically-modified 
twisted intertwining operators. 

We first have to recall briefly geometrically-modified twisted intertwining 
operators and pseudo-$q$-traces. 
Given a twisted intertwining operator $\Y$ of type $\binom{W_{3}}{W_{1}\;W_{2}}$
and $w_{1}\in W_{1}$, we have an operator (actually a series with linear maps from $W_{2}$
to $W_{3}$   as coefficients) $\mathcal{Y}_{1}(w_{1}, z)$.
The corresponding geometrically-modified operator is 
$$\mathcal{Y}_{1}(\mathcal{U}(q_{z})w_{1}, q_{z}),$$
 where 
$q^{z}=e^{2\pi iz}$,
$\mathcal{U}(q_{z})=(2\pi iq_{z})^{L(0)}e^{-L^{+}(A)}$ and 
$A_{j}\in \C$ for $j\in \mathbb{Z}_{+}$ are defined by 
\begin{eqnarray*}
\frac{1}{2\pi i}\log(1+2\pi i y)=\left(\exp\left(\sum_{j\in \mathbb{Z}_{+}}
A_{j}y^{j+1}\frac{\partial}{\partial y}\right)\right)y.
\end{eqnarray*}
See \cite{H-modular} for details. 

For pseudo-traces, we need to consider grading-restricted twisted $V$-modules 
equipped with a projective right module structure for an finite-dimensional associative algebra $P$ over $\C$. 
We first define pseudo-traces for a finitely generated projective right $P$-module $M$. 
For such a right $P$-module, there exists a projective basis, that is,
a pair of sets $\{w_i\}_{i=1}^n\subseteq M$,
$\{w'_i\}_{i=1}^n\subseteq \hom_P(M, P)$
such that for all $w\in M$, $w = \sum_{i=1}^n w_iw'_i(w).$
A linear function $\phi: P\rightarrow \C$ is said to be symmetric
if  $\phi(pq) = \phi(qp)$ for all $p,q\in P$.
For a symmetric linear function $\phi$,
the pseud-trace  $\tr^{\phi}_M\alpha$ for $\alpha\in {\rm End}_P(M)$ associated to $\phi$
is the function $\tr^{\phi}_M$ defined by
$$\tr^{\phi}_M \alpha = \phi\left(\sum_{i=1}^n w'_i(\alpha(w_i))\right).$$ 
For a grading-restricted twisted $V$-module $W$ 
equipped with a projective right $P$-module structure, its homogeneous subspaces $W_{[n]}$ 
for $n\in \C$ are finitely generated projective right $P$-modules. Then for a given 
symmetric linear function $\phi$ on $P$, we have the pseudo-trace $\tr^{\phi}_M\alpha_{n}$
of $\alpha_{n}\in {\rm End}_{P}(W_{[n]})$. For $\alpha\in {\rm End}_{P}(W)$,
we have $\alpha_{n}=\pi_{n}\alpha|_{W_{[n]}}\in {\rm End}_{P}(W_{[n]})$. We define 
$$\tr^{\phi}_W\alpha=\sum_{n\in \C}\tr^{\phi}_{W_{[n]}}\alpha_{n}.$$
Note that $\tr^{\phi}_W\alpha$ a a series of complex numbers, not a complex number. 
If we want to get a pseudo-trace in $\C$, we have to prove its convergence.

\vspace{1em}
\noindent {\bf Convergence and extension property of pseudo-$q$-traces of products of 
$n$ geomet\-rically-modified twisted intertwining operators}
Let $g_{i}$ for $i=1, \dots, n+1$ be automorphisms of $V$. Let 
$W_{i}$ for $i=1, \dots, n$ be grading-restricted generalized $g_{i}$-twisted $V$-modules, 
$\tilde{W}_{i}$ for $i=1, \dots, n$  $(g_{i+1}\cdots g_{n+1})$-twisted $V$-modules, and
$\mathcal{Y}_{i}$ for $i=1, \dots, n$ twisted intertwining operators of 
types $\binom{\tilde{W}_{i-1}}{ W_{i}\tilde{W}_{i}}$, respectively,
where we use the convention $\tilde{W}_{0}=
\tilde{W}_{n}$. 
Let $P$ be a finite-dimensional associative algebra and  $\phi$ a symmetric linear function
on $P$. Assume that $\tilde{W}_{0}=
\tilde{W}_{n}$ is also a projective right $P$-module and its twisted vertex operators
commute with the action of $P$. Assume in addition that the product 
$\mathcal{Y}_{1}(w_{1}, x_{1})
\cdots
\mathcal{Y}_{n}(w_{n}, x_{n})$ for $w_{1}\in W_{1}, \dots, w_{n}\in W_{n}$
commutes with the action of $P$. 
For $w_{i}\in W_{i}$, $i=1, \dots, n$,
$$\tr^{\phi}_{\tilde{W}_{n}}\mathcal{Y}_{1}(\mathcal{U}(q_{z_{1}})w_{1}, q_{z_{1}})
\cdots
\mathcal{Y}_{n}(\mathcal{U}(q_{z_{n}})w_{n}, q_{z_{n}})
q_{\tau}^{L(0)-\frac{c}{24}}$$
 is absolutely convergent in the region
  $1>|q_{z_1}| > \ldots > |q_{z_n}| > |q_\tau| > 0$ and can be extended
  to a multivalued analytic function 
$$
  \overline{F}^\phi_{\Y_1,\ldots, \Y_n}(w_1,\ldots, w_n; z_1,\ldots, z_n; \tau).
$$
in the region $\Im (\tau) >0$,
  $z_i \neq z_j + l + m\tau$ for $i\neq j$, $l,m\in \Z$.

\vspace{1em}

These  multivalued analytic functions are also conjectured to have associativity and  commutativity. 
These properties are consequences of the convergence and extension property 
 of pseudo-$q$-traces of products of 
$n$ geometrically-modified twisted intertwining operators and the  
associativity and  commutativity of twisted intertwining operators.

\vspace{1em}

\noindent {\bf Genus-one associativity}
In the setting of the convergence and extension property 
 of pseudo-$q$-traces of products of 
$n$ geometrically-modified twisted intertwining operators,  for $1\le k\le n-1$, 
there exist a $(g_{k}g_{k+1})$-twisted $V$-module $\hat{W}_{k}$ and twisted
intertwining operators $\hat{\mathcal{Y}}_{k}$ and 
$\hat{\mathcal{Y}}_{k+1}$ of 
types $\binom{\hat{W}_{k}}{W_{k} W_{k+1}}$ and 
$\binom{\tilde{W}_{k-1}}{\hat{W}_{k}\tilde{W}_{k+1}}$, respectively, 
such that 
\begin{align*}
\overline{F}^{\phi}&_{\mathcal{Y}_{1}, \dots, \mathcal{Y}_{k-1},
\hat{\mathcal{Y}}_{k+1}, \mathcal{Y}_{k+2}, \dots,
\mathcal{Y}_{n}}(w_{1}, 
\dots, w_{k-1}, \hat{\mathcal{Y}}(w_{k}, z_{k}-z_{k+1})
w_{k+1}, \nn
&\quad\quad\quad\quad\quad\quad\quad\quad\quad\quad\quad\quad\quad
w_{k+2}, \dots, w_{n};
z_{1}, \dots, z_{k-1}, z_{k+1}, \dots, z_{n}; \tau)
\end{align*}
is absolutely convergent in the region $1>|q_{z_{1}}|>\cdots 
>|q_{z_{k-1}}|>|q_{z_{k+1}}|>
\dots >|q_{z_{n}}|>|q_{\tau}|>0$ and 
$1>|q_{(z_{k}-z_{k+1})}-
1|>0$
and is convergent to 
$$\overline{F}^{\phi}_{\mathcal{Y}_{1}, \dots, \mathcal{Y}_{n}}(w_{1}, 
\dots, w_{n};
z_{1}, \dots, z_{n}; \tau)$$ 
in the region $1>|q_{ z_{1}}|>\cdots 
>|q_{z_{n}}|>|q_{\tau}|>0$ and 
$|q_{(z_{k}-z_{k+1})}|>1>|q_{(z_{k}-z_{k+1})}-
1|>0$.

\vspace{1em}

\noindent {\bf Genus-one commutativity}
In the setting of the convergence and extension property 
 of pseudo-$q$-traces of products of 
$n$ geometrically-modified twisted intertwining operators,  for $1\le k\le n-1$, there exist
a grading-restricted generalized $(g_{k}g_{k+2}\cdots g_{n+1})$-twisted  
$V$-module $\hat{W}_{k}$ and twisted intertwining operators
$\hat{\mathcal{Y}}_{k}$ and $\hat{\mathcal{Y}}_{k+1}$
of types $\binom{\hat{W}_{k}}{W_{k}\; \tilde{W}_{k+1}}$
and $\binom{\tilde{W}_{k-1}}{\phi_{g_{k}}(W_{k+1})\;\hat{W}_{k}}$, respectively,
or 
a grading-restricted generalized $(g_{k+1}^{-1}g_{k}g_{k+1}g_{k+2}\cdots g_{n+1})$-twisted  
$V$-module $\hat{W}_{k}$ and twisted intertwining operators
$\hat{\mathcal{Y}}_{k}$ and $\hat{\mathcal{Y}}_{k+1}$
of types $\binom{\hat{W}_{k}}{\phi_{g_{k+1}^{-1}}(W_{k})\; \tilde{W}_{k+1}}$
and $\binom{\tilde{W}_{k-1}}{W_{k+1}\;\hat{W}_{k}}$, respectively,
such that 
\begin{align*}
\overline{F}^{\phi}&_{\mathcal{Y}_{1}, \dots, \mathcal{Y}_{n}}(w_{1}, 
\dots, w_{n};
z_{1}, \dots, z_{n}; \tau)\nn
&=\overline{F}^{\phi}_{\mathcal{Y}_{1}, \dots, \mathcal{Y}_{k-1},
\hat{\mathcal{Y}}_{k+1}, \hat{\mathcal{Y}}_{k}, \mathcal{Y}_{k+2}
\dots, \mathcal{Y}_{n}}(w_{1}, 
\dots, w_{k-1}, w_{k+1}, w_{k}, w_{k+2}, \dots, w_{n};\nn
&\quad\quad\quad\quad\quad\quad\quad\quad\quad\quad\quad\quad
\quad\quad\quad
z_{1}, \dots, z_{k-1}, z_{k+1}, z_{k}, z_{k+2}, 
\dots, z_{n}; \tau).
\end{align*}

\vspace{1em}

The most important conjectural property of twisted intertwining operators
in the genus-one case is the modular invaraince. 

\vspace{1em}
\noindent {\bf Modular invariance of twisted intertwining operators}
For automorphisms $g_{i}$ of $V$, 
grading-restricted generalized $g_{i}$-twisted $V$-modules $W_{i}$ and $w_{i}\in W_{i}$
for $i=1, \dots, n$, let $\mathcal{F}_{w_{1}, \dots, w_{n}}$ be the vector space
spanned by functions of the form 
$$\overline{F}^\phi_{\Y_1,\ldots, \Y_n}(w_1,\ldots, w_n; z_1,\ldots, z_n; \tau)$$
for all finite-dimensional associative algebras $P$, all symmetric linear functions
$\phi$,  all $(g_{i+1}\cdots g_{n+1})$-twisted $V$-modules $\tilde{W}_{i}$ for $i=1, \dots, n+1$,
all projective right $P$-module structures on $\tilde{W}_{n+1}$ commuting with 
 its twisted vertex operators,
all twisted intertwining operators $\mathcal{Y}_{i}$ of 
types $\binom{\tilde{W}_{i-1}}{ W_{i}\tilde{W}_{i}}$  for $i=1, \dots, n$ , respectively,
such that their product commutes with the right action of $P$ on $\tilde{W}_{n+1}$. 
Then for 
$$\left(\begin{array}{cc}
\alpha&\beta\\
\gamma&\delta
\end{array}\right)\in SL(2, \mathbb{Z}),$$
$$
\overline{F}^{\phi}_{\mathcal{Y}_{1}, \dots, \mathcal{Y}_{n}}
\Biggl(\left(\frac{1}{\gamma\tau+\delta}\right)^{L(0)}w_{1}, \dots,
\left(\frac{1}{\gamma\tau+\delta}\right)^{L(0)}w_{n};
\frac{z_{1}}{\gamma\tau+\delta}, \dots, \frac{z_{n}}{\gamma\tau+\delta}; 
\frac{\alpha\tau+\beta}{\gamma\tau+\delta}\Biggr)$$
is in $\mathcal{F}_{w_{1}, \dots, w_{n}}$. 

\vspace{1em}

Since the commutativity, genus-one associativity and genus-one commutativity
are consequences of the other properties, the main properties that need to be proved
are the convergence and extension property of products of twisted intertwining operators,
associativity of twisted intertwining operators, convergence and extension property of pseudo-$q$-traces of 
products of $n$ geometrically-modified twisted intertwining operators and the modular invariance 
of twisted intertwining operators. 

\renewcommand{\theequation}{\thesection.\arabic{equation}}
\renewcommand{\thethm}{\thesection.\arabic{thm}}
\setcounter{equation}{0}
\setcounter{thm}{0}
\section{Some thoughts on further developments}

We discuss  in this section briefly some thoughts  of the author 
on further developments based on the conjectural properties  in the preceding section.

Assuming that 
the conjectural properties in Section 6 hold for twisted intertwining operators among suitable 
 twisted modules associated to a group of automorphisms  of a vertex operator algebra. 
Then we can generalize the results  described 
in Section 2 to results on what can be called ``equivariant chiral and full conformal field theories.''
In particular, we should be able to construct equivariant modular functors, equivariant genus-zero 
 and genus-one chrial conformal field theories. 
In Section 3, we mentioned that to construct higher-genus conformal field theories, the problem of
the convergence of multi pseudo-$q$-traces  of products of geometrically-modified 
intertwining operators is still open. 
In the case of orbifold conformal field theory, there is also a conjectural convergence of multi $q$-traces  
of products of geometrically-modified twisted  intertwining operators. If this convergence holds, then we can 
construct ``equivariant (all-genus) chiral and full conformal field theories."

We will also be able to obtain $G$-crossed tensor category structures. 
As is mentioned in Section 3, Conjecture \ref{g-crossed} is a consequence of Conjecture \ref{twisted-int-op}.
In fact, without assuming that Conjecture \ref{twisted-int-op} holds, 
tensor product bifunctors for suitable twisted modules  can be defined
in the same way as in \cite{tensor1}, \cite{tensor3} and
\cite{HLZ2} with intertwining operators replaced by twisted intertwining operators in \cite{H-twisted-int}. 
On the other hand, 
 though we do have a Jacobi identity for twisted vertex operators and one
twisted intertwining operator obtained as 
a special case for the Jacobi identity for intertwining operators in \cite{H-generalizedrationality}, it is not as nice 
as the one for vertex operators and one intertwining operator in \cite{FHL}. Thus
the construction of tensor product modules using the compatibility condition
and local grading-restriction condition has to be modified by using the method of formulating and studying 
twisted intertwining operators in \cite{H-twisted-int} instead of the method based on the Jacobi identity 
in \cite{tensor1}, \cite{tensor2}, \cite{tensor3} and \cite{HLZ3}. 
Now assuming that  the convergence and extension property of products of twisted intertwining operators
and associativity of twisted intertwining operators hold. Then the associativity of 
the tensor product bifunctors can be proved in the same way as in \cite{tensor4} and 
\cite{HLZ6}. With this construction, we actually obtain what should be called 
a ``vertex $G$-crossed tensor category'' structure and the 
$G$-crossed tensor category structure can be derived from this structure in the same way as 
how the braided tensor category structure is derived from the vertex tensor category structure 
in \cite{H-rigidity} and \cite{HLZ8}. We also expect that if 
the convergence and  extension property of  $q$-traces of 
products of $n$ geometrically-modified twisted intertwining operators and 
the modular invariance of twisted intertwining operators hold, then the rigidity 
can be proved in the same way as in \cite{H-rigidity}  in the case that the category 
of lower-bounded generalized twisted modules are semisimple.

The conjectural properties in Section 6 can be proved if we know that the fixed point subalgebra $V^{G}$ 
of $V$ satisfies conditions in the papers \cite{H-diff-eqn}, \cite{H-modular}, \cite{H-verlinde-conj}
and \cite{H-rigidity} or some other conditions. 
See \cite{Mi} and \cite{CM} for proofs of the conditions for $V^{G}$ in the case that 
$G$ is finite cyclic and \cite{Mc1} and \cite{Mc2} for some results one can obtain when 
$V^{G}$ is assumed to satisfy certain conditions, including in particular suitable conditions 
on suitable categories of $V^{G}$-modules. But proving these conditions to hold for $V^{G}$  seems to be at least as
difficult as proving the conjectural properties  in Section 6. Also note that the reason why we want to 
prove these conditions for $V^{G}$ is exactly that these conditions can be used to prove the
the conjectural properties  in Section 6. In fact, the proof of the reductivity property of $V^{G}$ 
in \cite{CM} uses heavily the theory of intertwining operators established in \cite{H-diff-eqn}
and \cite{H-modular}. The author believes that these conditions for $V^{G}$  follow from 
the conjectures in Section 3. Therefore in the author's opinion, if possible, we should prove
the conjectural properties  in Section 6 directly, without using the fixed point subalgebra $V^{G}$. 
Then we should be able to derive the algebraic  conditions for $V^{G}$ as consequences. In order to derive the
conditions on $V^{G}$ as consequences, we need to prove that for a general vertex operator algebra $V$,
the properties of (untwisted)  intertwining operators (such as the associativity and modular invariance) 
imply the algebraic conditions on $V$. This problem is  interesting
even without the study of orbifold conformal field theory since it provides a deep understanding 
of the connection between algebraic properties of the vertex operator algebra $V$ and 
geometric properties of genus-zero and genus-one correlation functions of the corresponding conformal 
field theory.

In \cite{H-cft-lattice}, the author introduced a precise notion of dual of an intertwining operator algebra. 
Since vertex operator algebra is also an intertwining operator algebra, we also have a dual of 
a vertex operator algebra. In fact, the dual of a vertex operator algebra is simply 
the intertwining operator algebra constructed using all intertwining operators 
among all modules for the vertex operator algebras. In particular, a self-dual vertex operator algebra
in this sense means that the only irreducible module is the vertex operator algebra itself and 
is called in many papers a holomorphic vertex operator algebra. The moonshine module constructed in 
\cite{FLM} is a self-dual vertex operator algebra and the uniqueness conjecture of the moonshine
module states that a self-dual vertex operator algebra of central charge $24$ and without nonzero weight $1$ elements
must be isomorphic to the moonshine module as a vertex operator algebra. 
Let $V$ be a self-dual vertex operator algebra of central charge $24$ and without nonzero weight $1$ elements.
It is desirable if we can embed $V$ into a largest possible intertwining operator algebra.
Then we can perform all types of operations and constructions in this largest intertwining 
operator algebra containing $V$. But taking the dual of $V$ does not work since 
it is self dual. On the other hand, if we let $G$ be the full automorphism group of $V$, then we can take the dual of 
the vertex operator subalgebra $V^{G}$ of $V$. This dual of $V^{G}$ must contain $V$ and is 
the largest among the duals of all the fixed point subalgebras of $V$. The construction of the dual of $V^{G}$ is 
in fact part of Problem \ref{main-prob}. As is mentioned after Problem \ref{main-prob} and in the beginning of 
Section 5,  instead of studying directly $V^{G}$-modules and intertwining 
operators among $V^{G}$-modules, we study twisted $V$-modules and twisted intertwining operators
among twisted $V$-modules. Thus the construction and study of the dual of $V^{G}$ 
are reduced to the construction and study of the orbifold conformal field theory obtained from 
the vertex operator algebra $V$ and its automorphsim group $G$. We hope that 
if this orbifold conformal field theory is constructed, we can use it to find a strategy 
to prove the uniqueness of the moonshine module. Note that the construction of 
such an orbifold conformal field theory from such an arbitrary vertex operator algebra satisfying 
three conditions must be general and abstract and cannot be worked out for a particular example
such as the moonshine module itself. This is in fact one important reason why we have to 
develop a general orbifold conformal field theory instead of just explicit examples.

\noindent {\small \sc Department of Mathematics, Rutgers University,
110 Frelinghuysen Rd., Piscataway, NJ 08854-8019}

\noindent {\em E-mail address}: yzhuang@math.rutgers.edu

\end{document}